\documentclass[a4paper]{article}
\usepackage{amssymb}
\usepackage{amsfonts}
\usepackage{amsmath}

\input{tcilatex}

\begin{document}

\title{Schubert claculus and the Hopf algebra structures of exceptional Lie
groups}
\author{Haibao Duan and Xuezhi Zhao \\
Institute of Mathematics, Chinese Academy of Sciences,\\
Department of Mathematics, Capital Normal University}
\maketitle

\begin{abstract}
Let $G$ be an exceptional Lie group with a maximal torus $T$. Based on
common properties in the Schubert presentation of the cohomology ring $%
H^{\ast }(G/T;\mathbb{F}_{p})$ \cite{DZ1}, and concrete expressions of
generalized Weyl invariants for $G$ over $\mathbb{F}_{p}$, we obtain a
unified approach to the structure of $H^{\ast }(G;\mathbb{F}_{p})$ as a Hopf
algebra over the Steenrod algebra $\mathcal{A}_{p}$.

The results has been applied in \cite{Du2} to determine the near--Hopf ring
structure on the integral cohomology of all exceptional Lie groups.

\begin{description}
\item \textsl{2000 Mathematical Subject Classification:} 14M15; 57T15;

\item \textsl{Key words and phrases: }Lie groups; Schubert classes; Hopf
algebra; Steenrod algebra

\item \textsl{Email addresses:} dhb@math.ac.cn; zhaoxve@mail.cnu.edu.cn
\end{description}
\end{abstract}

\section{Introduction}

Let $G$ be an exceptional Lie group with a maximal torus $T\subset G$, a set
of fundamental dominant weights $\{\omega _{1},\cdots ,\omega _{n}\}\in
H^{2}(G/T)$ \cite{BH}, $n=\dim T$. Based on Schubert presentation for the
integral cohomology ring $H^{\ast }(G/T)$ of the flag manifold $G/T$
obtained in \cite{DZ1}, we have derived for each exceptional $G$ and prime $%
p=2,3,5$ an explicitly set $\{\theta _{s_{1}},\cdots ,\theta _{s_{n}}\}$, $%
\deg \theta _{s}=s$, of \textsl{generalized Weyl invariants} \textsl{over} $%
\mathbb{F}_{p}$ (see \S 2 and \S 5.2). In this paper we show how this set of
polynomials in $\omega _{1},\cdots ,\omega _{n}$ gives rise naturally to a
set $\{\alpha _{2s_{1}-1},\cdots ,\alpha _{2s_{n}-1}\}$ of $p$--\textsl{%
transgressive generators} on $H^{\ast }(G;\mathbb{F}_{p})$, and the
structure of $H^{\ast }(G;\mathbb{F}_{p})$ as a Hopf algebra over the
Steenrod algebra can be effectively calculated by computing with these
polynomials. We shall restrict ourself to the cases where the integral
cohomology $H^{\ast }(G)$ contains non--trivial $p$--torsion subgroup, for
exactly in these cases the Hopf algebra $H^{\ast }(G;\mathbb{F}_{p})$ fails
to be a primitive generated exterior algebra.

Let $\mathcal{A}_{{p}}$ be the $\func{mod}$ $p$ \textsl{Steenrod algebra}
with $\mathcal{P}^{k}\in \mathcal{A}_{{p}}$, $k\geq 1$, the $k^{th}$ \textsl{%
reduced power} \cite{SE} and $\delta _{p}\in \mathcal{A}_{{p}}$ the \textsl{%
Bockstein\ operator}. If $p=2$ it is also customary to write $Sq^{2k}$
instead of $\mathcal{P}^{k}$, $Sq^{1}$ in the place of $\delta _{2}$. Our
first result gives a complete characterization of $H^{\ast }(G;\mathbb{F}%
_{p})$ as an algebra over $\mathcal{A}_{{p}}$ in term of the $p$%
--transgressive generators $\{\alpha _{2s_{1}-1},\cdots ,\alpha
_{2s_{n}-1}\} $ on $H^{\ast }(G;F_{p})$.

\noindent \textbf{Theorem 1.} \textsl{Let }$(G,p)$\textsl{\ be a pair with }$%
G$\textsl{\ an exceptional Lie group and }$H^{\ast }(G)$\textsl{\ containing
non--trivial }$p$\textsl{--torsion subgroup..}

\begin{enumerate}
\item[(1.1)] \textsl{The algebra }$H^{\ast }(G;\mathbb{F}_{2})$ \textsl{has
the presentation}
\end{enumerate}

$H^{\ast }(G_{2};\mathbb{F}_{2})=\mathbb{F}_{2}[x_{6}]/\left\langle
x_{6}^{2}\right\rangle \otimes \Delta _{\mathbb{F}_{2}}(\alpha _{3},\alpha
_{5})$;

$H^{\ast }(F_{4};\mathbb{F}_{2})=\mathbb{F}_{2}[x_{6}]/\left\langle
x_{6}^{2}\right\rangle \otimes \Delta _{\mathbb{F}_{2}}(\alpha _{3},\alpha
_{5},\alpha _{15},\alpha _{23})$;

$H^{\ast }(E_{6};\mathbb{F}_{2})=\mathbb{F}_{2}[x_{6}]/\left\langle
x_{6}^{2}\right\rangle \otimes \Delta _{\mathbb{F}_{2}}(\alpha _{3},\alpha
_{5},\alpha _{9},\alpha _{15},\alpha _{17},\alpha _{23})$;

$H^{\ast }(E_{7};\mathbb{F}_{2})=\frac{\mathbb{F}_{2}[x_{6},x_{10},x_{18}]}{%
\left\langle x_{6}^{2},x_{10}^{2},x_{18}^{2}\right\rangle }\otimes \Delta _{%
\mathbb{F}_{2}}(\alpha _{3},\alpha _{5},\alpha _{9},\alpha _{15},\alpha
_{17},\alpha _{23},\alpha _{27})$;

$H^{\ast }(E_{8};\mathbb{F}_{2})=\frac{\mathbb{F}%
_{2}[x_{6},x_{10},x_{18},x_{30}]}{\left\langle
x_{6}^{8},x_{10}^{4},x_{18}^{2},x_{30}^{2}\right\rangle }\otimes \Delta _{%
\mathbb{F}_{2}}(\alpha _{3},\alpha _{5},\alpha _{9},\alpha _{15},\alpha
_{17},\alpha _{23},\alpha _{27},\alpha _{29})$,

\noindent \textsl{on which} \textsl{all nontrivial }$\mathcal{A}_{2}$\textsl{%
--actions} \textsl{are given by}

\begin{quote}
$\delta _{2}(\alpha _{5})=x_{6}$ \textsl{in} $G_{2},F_{4},E_{6},E_{7},E_{8}$;

$\delta _{2}(\alpha _{2r-1})=x_{2r},$ $r=5,9$;

$\delta _{2}(\alpha _{15})=x_{6}x_{10}$; $\quad \delta _{2}(\alpha
_{27})=x_{10}x_{18}$ \textsl{in} $E_{7},E_{8}$

$\delta _{2}(\alpha _{23})=x_{6}x_{18}$ \textsl{in} $E_{7}$;

$\delta _{2}(\alpha _{23})=x_{6}x_{18}+x_{6}^{4}$; $\quad \delta _{2}(\alpha
_{29})=x_{30}+x_{6}^{2}x_{18}$ \textsl{in} $E_{8}$,

$\mathcal{P}^{1}\alpha _{3}=\alpha _{5}$ \textsl{in} $%
G_{2},F_{4},E_{6},E_{7},E_{8}$;

$\mathcal{P}^{4}\alpha _{15}=\alpha _{23}$ \textsl{in} $%
F_{4},E_{6},E_{7},E_{8}$;

$\mathcal{P}^{2}\alpha _{5}=\alpha _{9}$; $\mathcal{P}^{4}\alpha _{9}=%
\mathcal{P}^{1}\alpha _{15}=\alpha _{17}$ \textsl{in} $E_{6},E_{7},E_{8}$;

$\mathcal{P}^{2}\alpha _{23}=\alpha _{27}$ \textsl{in} $E_{7},E_{8}$;

$\mathcal{P}^{3}\alpha _{23}=\mathcal{P}^{1}\alpha _{27}=\alpha _{29}$ 
\textsl{in} $E_{8}$;
\end{quote}

\begin{enumerate}
\item[(1.2)] \textsl{The algebra }$H^{\ast }(G;\mathbb{F}_{3})$ \textsl{has
the presentation}
\end{enumerate}

\begin{quote}
$H^{\ast }(F_{4};\mathbb{F}_{3})=\mathbb{F}_{3}[x_{8}]/\left\langle
x_{8}^{3}\right\rangle \otimes \Lambda _{\mathbb{F}_{3}}(\alpha _{3},\alpha
_{7},\alpha _{11},\alpha _{15})$;

$H^{\ast }(E_{6};\mathbb{F}_{3})=\mathbb{F}_{3}[x_{8}]/\left\langle
x_{8}^{3}\right\rangle \otimes \Lambda _{\mathbb{F}_{3}}(\alpha _{3},\alpha
_{7},\alpha _{9},\alpha _{11},\alpha _{15},\alpha _{17})$;

$H^{\ast }(E_{7};\mathbb{F}_{3})=\mathbb{F}_{3}[x_{8}]/\left\langle
x_{8}^{3}\right\rangle \otimes \Lambda _{\mathbb{F}_{3}}(\alpha _{3},\alpha
_{7},\alpha _{11},\alpha _{15},\alpha _{19},\alpha _{27},\alpha _{35})$;

$H^{\ast }(E_{8};\mathbb{F}_{3})=\mathbb{F}_{3}[x_{8},x_{20}]/\left\langle
x_{8}^{3},x_{20}^{3}\right\rangle \otimes \Lambda _{\mathbb{F}_{3}}(\alpha
_{3},\alpha _{7},\alpha _{15},\alpha _{19},\alpha _{27},\alpha _{35},\alpha
_{39},\alpha _{47})$
\end{quote}

\noindent \textsl{on which} \textsl{all nontrivial }$\mathcal{A}_{3}$\textsl{%
--actions} \textsl{are given by}

\begin{quote}
$\delta _{3}(\alpha _{7})=-x_{8}$; $\quad \delta _{3}(\alpha
_{15})=-x_{8}^{2}$ \textsl{in} $F_{4},E_{6},E_{7},E_{8}$;

$\delta _{3}(\alpha _{19})=x_{20}$; $\quad \delta _{3}(\alpha
_{27})=-x_{8}x_{20}$; $\quad \delta _{3}(\alpha _{35})=x_{8}^{2}x_{20}$;

$\delta _{3}(\alpha _{39})=x_{20}^{2}$; $\quad \delta _{3}(\alpha
_{47})=x_{8}x_{20}^{2}$ \textsl{in} $E_{8}$

$\mathcal{P}^{1}\alpha _{3}=\alpha _{7}$ \textsl{in} $%
F_{4},E_{6},E_{7},E_{8} $;

$\mathcal{P}^{1}\alpha _{11}=\alpha _{15}$ \textsl{in} $F_{4},E_{6},E_{7}$;

$\mathcal{P}^{1}\alpha _{15}=\mathcal{P}^{3}\alpha _{7}=\alpha _{19}$; $%
\mathcal{P}^{3}\alpha _{15}=\alpha _{27}$ \textsl{in} $E_{7},E_{8}$;

$\mathcal{P}^{2}\alpha _{11}=-\alpha _{19}$ \textsl{in} $E_{7}$;

$\mathcal{P}^{1}\alpha _{35}=\mathcal{P}^{3}\alpha _{27}=\alpha _{39}$; $%
\mathcal{P}^{3}\alpha _{35}=\alpha _{47}$ \textsl{in} $E_{8}$;
\end{quote}

\begin{enumerate}
\item[(1.3)] \textsl{The algebra }$H^{\ast }(E_{8};\mathbb{F}_{5})$ \textsl{%
has the presentation}
\end{enumerate}

\begin{quote}
$\mathbb{F}_{5}[x_{12}]/\left\langle x_{12}^{5}\right\rangle \otimes \Lambda
(\alpha _{3},\alpha _{11},\alpha _{15},\alpha _{23},\alpha _{27},\alpha
_{35},\alpha _{39},\alpha _{47})$
\end{quote}

\noindent \textsl{on which} \textsl{all nontrivial }$\mathcal{A}_{5}$\textsl{%
--actions} \textsl{are given by}

\begin{quote}
$\delta _{5}(\alpha _{11})=-x_{12}$; $\delta _{5}(\alpha _{23})=-x_{12}^{2}$%
; $\delta _{5}(\alpha _{35})=x_{12}^{3}$; $\delta _{5}(\alpha
_{47})=2x_{12}^{4}$

$\mathcal{P}^{1}\alpha _{i}=\alpha _{i+8}$, $i=3,15,27,39$.
\end{quote}

\noindent \textbf{Remark 1.1.} In the classical descriptions of $H^{\ast
}(E_{7};\mathbb{F}_{2})$ and $H^{\ast }(E_{8};\mathbb{F}_{2})$ in \cite{Ar1,
AS, T, Ka, Ko1, KN, Mi} the generators were specified mainly up to their
degrees and the action of $Sq^{1}=\delta _{2}$ on the generators in degrees $%
15,23,27$ was absent. In comparison, results in (1.1) constitutes a complete
characterization of $H^{\ast }(G;\mathbb{F}_{2})$ as an algebra over $%
\mathcal{A}_{2}$, see Corollary 4.2 and Remark 4.3.

In \cite{KM3} Kono and Mimura largely determined the $\mathcal{A}_{3}$
action on $H^{\ast }(E_{7};\mathbb{F}_{3})$ and $H^{\ast }(E_{8};\mathbb{F}%
_{3})$ with respect also to a set of transgressive generators, except an
indeterminacy $\epsilon =\pm 1$ occurred in their expressions of $\mathcal{P}%
^{2}e_{11},\mathcal{P}^{1}e_{15}$ in $E_{7}$, and of $\mathcal{P}^{1}e_{15}$%
, $\mathcal{P}^{1}e_{35}$ in $E_{8}$. Again, with respect to our explicit
construction these ambiguities are clarified in (1.2).

Result (1.3) essentially agrees with the calculation of Kono \cite[Theorem
5.15]{Ko2}, whose generators $x_{3},x_{15},x_{27},x_{39}$ are also
transgressive, and correspond to ours $\alpha _{3},2\alpha _{15},3\alpha
_{27},2\alpha _{39}$, respectively.$\square $

\bigskip

For a prime $p$ the multiplication $\mu :G\times G\rightarrow G$ in $G$
induces the \textsl{reduced} \textsl{co--product}

\begin{center}
$\phi _{p}:H^{\ast }(G;\mathbb{F}_{p})\rightarrow H^{\ast }(G;\mathbb{F}%
_{p})\otimes H^{\ast }(G;\mathbb{F}_{p})$
\end{center}

\noindent by $\phi _{p}(x)=\mu ^{\ast }(x)-(x\otimes 1+1\otimes x)$ by
virtue of the Kunneth formula. It furnishes $H^{\ast }(G;\mathbb{F}_{p})$
with the structure of a \textsl{Hopf algebra}. The proof of the next result
reveals that, with respect to our presentation of $H^{\ast }(G;\mathbb{F}%
_{p})$ in Theorem 1, \textsl{the structure of }$H^{\ast }(G;\mathbb{F}_{p})$%
\textsl{\ as a Hopf algebra is entirely determined by its structure as an
algebra over }$\mathcal{A}_{p}$.

\bigskip

\noindent \textbf{Theorem 2. }\textsl{With respect to the presentation of }$%
H^{\ast }(G;\mathbb{F}_{p})$\textsl{\ in (1.1)--(1.3), the reduced
co--product }$\phi _{p}$\textsl{\ on }$H^{\ast }(G;\mathbb{F}_{p})$\textsl{\
is determined, respectively, by}

\begin{enumerate}
\item[(1.4)] \textsl{for all exceptional }$G$\textsl{\ and }$p=2:$

$\phi _{2}(\alpha _{i})=\left\{ 
\begin{tabular}{l}
$0$ \textsl{for} $i=3$, \textsl{or} $i=15$ \textsl{and} $G=F_{4}$; \\ 
$x_{6}\otimes \alpha _{9}$ \textsl{if} $i=15$ \textsl{and} $G=E_{6}$; \\ 
$x_{10}\otimes \alpha _{5}+x_{6}\otimes \alpha _{9}$ \textsl{if }$i=15$%
\textsl{\ and} $G=E_{7}$; \\ 
$x_{10}\otimes \alpha _{5}+x_{6}\otimes \alpha _{9}+x_{6}^{2}\otimes \alpha
_{3}$ \textsl{if }$i=15$\textsl{\ and} $G=E_{8}$,%
\end{tabular}%
\ \ \right. $

\item[(1.5)] \textsl{for an exceptional }$G\neq G_{2}$\textsl{\ and }$p=3:$

$\phi _{3}(\alpha _{i})=\left\{ 
\begin{tabular}{l}
$0$ \textsl{if} $i=3,7,9,17,19$ \\ 
$-x_{8}\otimes \alpha _{3}$ \textsl{if }$i=11$\textsl{\ and} $%
G=F_{4},E_{6},E_{7}$ \\ 
$-x_{8}\otimes \alpha _{7}$ \textsl{if} $(G,i)=(E_{8},15)$; \\ 
$x_{8}\otimes \alpha _{27}+x_{8}^{2}\otimes \alpha _{19}$ \textsl{if} $%
(G,i)=(E_{7},35)$; \\ 
${\small x}_{8}{\small \otimes \alpha }_{27}{\small +x}_{8}^{2}{\small %
\otimes \alpha }_{19}{\small +x}_{8}{\small x}_{20}{\small \otimes \alpha }%
_{7}{\small -x}_{20}{\small \otimes \alpha }_{15}$ \textsl{if} $%
(G,i)=(E_{8},35)$,%
\end{tabular}%
\ \ \ \right. $

\item[(1.6)] \textsl{for }$(G,p)=(E_{8},5):$
\end{enumerate}

\begin{center}
$\phi _{5}(\alpha _{i})=\left\{ 
\begin{tabular}{l}
$0$ \textsl{if} $i=3$, \\ 
$2x_{12}\otimes \alpha _{3}$ \textsl{if} $i=15$, \\ 
$2x_{12}\otimes \alpha _{15}+2x_{12}^{2}\otimes \alpha _{3}$ \textsl{if} $%
i=27$, \\ 
$3x_{12}\otimes \alpha _{27}+3x_{12}^{2}\otimes \alpha
_{15}+2x_{12}^{3}\otimes \alpha _{3}$ \textsl{if} $i=39$.%
\end{tabular}%
\ \ \ \ \ \right. $
\end{center}

Historically, the Hopf algebras $H^{\ast }(G;\mathbb{F}_{p})$ for
exceptional Lie groups $G$ were studied by quite different methods,
presented by generators with various origins and using case by case
computations depending on $G$ and $p$. As examples, see

\qquad Borel \cite{B} for $(G,p)=(G_{2};2),(F_{4};2)$;

\qquad Araki \cite{Ar2} for $(F_{4};3)$;

\qquad Toda, Kono, Mimura, Shimada \cite{To,KM1, KMS} for $(E_{i};2)$, $%
i=6,7,8$;

\qquad Kono, Mimura and Toda \cite{KM2, To} for $(E_{6};3)$;

\qquad Kono-Mimura \cite{KM3} for $(E_{7};3)$ and $(E_{8};3)$;

\qquad Kono \cite{Ko2} for $(E_{8};5)$.

\noindent "\textsl{All these works created an impression that there is no
simple procedure for calculating} $H^{\ast }(G;\mathbb{F}_{p})$", as
commented by Ka\v{c} \cite{K}. In Ka\v{c} \cite{K} initiated a unified
approach to $H^{\ast }(G;\mathbb{F}_{p})$ in the context of Schubert
calculus on $G/T$. He showed that the algebra structure for $p\neq 2$ and
additive structure for $p=2$ are entirely determined by the degrees of basic
Weyl--invariants over $\mathbb{Q}$, and the degrees $\{\deg \theta
_{s_{i}}\mid 1\leq i\leq n\}$ of the basic \textquotedblleft generalized
invariants\textquotedblright\ over $\mathbb{F}_{p}$. Indeed, with concrete
expressions for the set of polynomials $\{\theta _{s_{i}}\mid 1\leq i\leq
n\} $ we have arrived at the structure of $H^{\ast }(G;\mathbb{F}_{p})$ as a
Hopf algebra over $\mathcal{A}_{p}$.

In the modern form taken by the classical topic known as "\textsl{%
enumerative geometry}" in the 19$^{th}$ century, Schubert calculus amounts
to calculation in the intersection ring $H^{\ast }(G/T)$ of the flag
manifold $G/T$ \cite[p.331]{W}. The present work, as well as the relevant
ones \cite{DZ1, DZ2, DZ3, Du2}, serves the purpose to demonstrate that the
cohomology theory of Lie groups can be boiled down to computation with
certain polynomials in Schubert classes.

This paper is arranged as follows. In \S 2 we recall from \cite{DZ1} common
properties in the Schubert presentation of the ring $H^{\ast }(G/T;\mathbb{F}%
_{p})$. In particular, the role of the set $\{\theta _{s_{1}},\cdots ,\theta
_{s_{n}}\}$ of generalized Weyl invariants is emphasized. In \S 3 we
construct $H^{\ast }(G;\mathbb{F}_{p})$ from the presentation of $H^{\ast
}(G/T;\mathbb{F}_{p})$ in \S 2, and relate the $\mathcal{P}^{k}$ action on $%
H^{\ast }(G;\mathbb{F}_{p})$ to certain relations among the polynomials $%
\theta _{s_{1}},\cdots ,\theta _{s_{n}}$. \S 5 is created to record
expressions of these polynomials and to handle computational aspects in this
paper. With these preliminaries Theorems 1 and 2 are established in \S 4.

\section{Schubert presentation of $H^{\ast }(G/T;\mathbb{F}_{p})$}

For a Lie group $G$ with a maximal torus $T$ consider the fibration

\begin{enumerate}
\item[(2.1)] $G/T\overset{\psi }{\hookrightarrow }BT\overset{\pi }{%
\rightarrow }BG$
\end{enumerate}

\noindent induced by the inclusion $T\subset G$, where $BT$ (resp. $BG$) is
the classifying space of $T$ (resp. $G$). Since $H^{odd}(BT)$ $%
=H^{odd}(G/T)=0$, the cohomology exact sequence of the pair $(BT,G/T)$ in
the $\mathbb{F}_{p}$--coefficients contains the section

\begin{enumerate}
\item[(2.2)] $0\rightarrow H^{even}(BT,G/T;\mathbb{F}_{p})\overset{j}{%
\rightarrow }H^{\ast }(BT;\mathbb{F}_{p})\overset{\psi _{p}^{\ast }}{%
\rightarrow }H^{\ast }(G/T;\mathbb{F}_{p})$
\end{enumerate}

\noindent where, as is classical, $H^{\ast }(BT;\mathbb{F}_{p})$ can be
identified with the free polynomial ring $\mathbb{F}_{p}[\omega _{1},\cdots
,\omega _{n}]$ in a set of fundamental dominant weights $\omega _{1},\cdots
,\omega _{n}\in H^{2}(BT;\mathbb{F}_{p})$ of $G$, and where the\textbf{\ }%
ring map $\psi _{p}^{\ast }$ induced by the fiber inclusion $\psi $ is well
known as the \textsl{Borel's characteristic map} in characteristic $p$ \cite%
{BH, B1}.

According to Borel \cite{B1}, if the integral cohomology $H^{\ast }(G)$ is
free of $p$--torsion, then $\psi _{p}^{\ast }$ is surjective and induces an
isomorphism

\begin{quote}
$H^{\ast }(G/T;\mathbb{F}_{p})=H^{\ast }(BT;\mathbb{F}_{p})/\left\langle
H^{+}(BT;\mathbb{F}_{p})^{W(G)}\right\rangle $
\end{quote}

\noindent where $\left\langle H^{+}(BT;\mathbb{F}_{p})^{W(G)}\right\rangle $
is the ideal in $H^{\ast }(BT;\mathbb{F}_{p})$ generated by Weyl invariants
in positive degrees (see Demazure \cite{D} for another proof of this fact).
Without any restriction on the torsion subgroup of $H^{\ast }(G)$ we extends
this classical result in Lemma 2.1 below.

For simplicity, we make no difference in notation between a polynomial $%
\theta \in H^{\ast }(BT;\mathbb{F}_{p})$ and its $\psi _{p}^{\ast }$--image
in $H^{\ast }(G/T;\mathbb{F}_{p})$. Given a subset $\{f_{1},\ldots ,f_{m}\}$
in a ring write $\left\langle f_{1},\ldots ,f_{m}\right\rangle $ for the
ideal generated by $f_{1},\ldots ,f_{m}$.

\bigskip

\noindent \textbf{Lemma 2.1. }(\cite[Proposition 3]{DZ1})\textbf{.} \textsl{%
For each }$1$\textsl{--connected Lie group }$G$\textsl{\ with rank }$n$%
\textsl{\ and and a prime }$p$\textsl{, there exist}

\begin{quote}
\textsl{a set }$\{\theta _{s_{1}},\cdots ,\theta _{s_{n}}\}\subset H^{\ast
}(BT;\mathbb{F}_{p})$\textsl{\ of }$n$\textsl{\ polynomials; and}

\textsl{a set }$\{y_{t_{1}},\cdots ,y_{t_{k}}\}\subset H^{\ast }(G/T;\mathbb{%
F}_{p})$\textsl{\ of Schubert classes on }$G/T$
\end{quote}

\noindent \textsl{with }$\deg \theta _{s}=2s$\textsl{, }$\deg y_{t}=2t>2$%
\textsl{, so that}

\textsl{i) }$\ker \psi _{p}^{\ast }=\left\langle \theta _{s_{1}},\cdots
,\theta _{s_{n}}\right\rangle $\textsl{;}

\textsl{ii) }$H^{\ast }(G/T;\mathbb{F}_{p})=\mathbb{F}_{p}[\omega
_{1},\ldots ,\omega _{n},y_{t}]/\left\langle \theta _{s},y_{t}^{k_{t}}+\beta
_{t}\right\rangle _{s\in r(G,p),\text{ }t\in e(G,p)}$ \textsl{with}

\noindent $\beta _{t}\in \left\langle \omega _{1},\ldots ,\omega
_{n}\right\rangle $\textsl{;}

\textsl{iii) the three sets\ of integers}

\begin{quote}
$r(G,p)=\{s_{1},\cdots ,s_{n}\}$\textsl{, }$e(G,p)=\{t_{1},\cdots ,t_{k}\}$%
\textsl{\ and }$\{k_{t}\}_{t\in e(G,p)}$
\end{quote}

\noindent \textsl{are subject to the constraints}

\begin{quote}
$e(G,p)\subset r(G,p)$\textsl{;\qquad }$\dim G=\tsum\limits_{s\in
r(G,p)}(2s-1)+\tsum\limits_{t\in e(G,p)}2(k_{t}-1)t$\textsl{.}$\square $
\end{quote}

Since the set $\left\{ \omega _{1},\cdots ,\omega _{n}\right\} $ of weights
consists of all Schubert classes on $G/T$ with cohomology degree $2$, ii) of
Lemma 2.1 describes the ring $H^{\ast }(G/T;\mathbb{F}_{p})$ by certain
Schubert classes on $G/T$ and therefore, is called \textsl{a Schubert
presentation} of $H^{\ast }(G/T;\mathbb{F}_{p})$. In addition to $\left\{
\omega _{1},\cdots ,\omega _{n}\right\} $ elements in the set $%
\{y_{t}\}_{t\in e(G,p)}$ will be called the $p$\textsl{--special Schubert
classes} on $G/T$. For each exceptional $G$ and prime $p$, a set of $p$%
\textsl{--}special Schubert classes on $G/T$ has been determined in \cite%
{DZ1}, and is specified by their Weyl coordinates in the table below:

\begin{center}
\begin{tabular}{|l|l|l|l|l|}
\hline
$y_{i}$ & $G_{2}/T$ & $F_{4}/T$ & $E_{n}/T,\text{ }n=6,7,8$ & $p$ \\ 
\hline\hline
$y_{3}$ & $\sigma _{\lbrack 1,2,1]}$ & $\sigma _{\lbrack 3,2,1]}$ & $\sigma
_{\lbrack 5,4,2]}\text{, }n=6,7,8$ & $2$ \\ 
$y_{4}$ &  & $\sigma _{\lbrack 4,3,2,1]}$ & $\sigma _{\lbrack 6,5,4,2]}\text{%
, }n=6,7,8$ & $3$ \\ 
$y_{5}$ &  &  & $\sigma _{\lbrack {7,6,5,4,2}]}\text{, }n=7,8$ & $2$ \\ 
$y_{6}$ &  &  & $\sigma _{\lbrack 1,3,6,5,4,2]}\text{, }n=8$ & $5$ \\ 
$y_{9}$ &  &  & $\sigma _{\lbrack 1,5,4,3,7,6,5,4,2]}$,$\text{ }n=7,8$ & $2$
\\ 
$y_{10}$ &  &  & $\sigma _{\lbrack {1,6,5,4,3,7,6,5,4,2}]}\text{, }n=8$ & $3$
\\ 
$y_{15}$ &  &  & $\sigma _{\lbrack 5,4,2,3,1,6,5,4,3,8,7,6,5,4,2]}\text{, }%
n=8$ & $2$ \\ \hline
\end{tabular}

{\small The }${\small p}$--{\small special Schubert classes on }${\small G/T}
${\small \ and their abbreviations}
\end{center}

In view of i) of Lemma 2.1 we shall call $\{\theta _{s}\}_{s\in r(G,p)}$ a
set of \textsl{generating polynomials} for $\ker \psi _{p}^{\ast }$. These
polynomials have been emphasized by Ka\v{c} \cite{K} as \textsl{a regular
sequence of homogeneous generators} for $\ker \psi _{p}^{\ast }$ (or a set
of \textsl{generalized Weyl invariants} for $G$ over $\mathbb{F}_{p}$);
notified by Ishitoya, Kono and Toda \cite{IKT} as the \textsl{transgressive
imagines} of a set of transgressive generators on $H^{\ast }(G;\mathbb{F}%
_{p})$. However, it is in the context of \cite[\S 6]{DZ1} that concrete%
\textbf{\ }expression of a set of\textsl{\ }such polynomials is available
for each exceptional $G$ and prime $p$.

Assume in the remainder of this section that $(G,p)$ is a pair with $G$
exceptional and $H^{\ast }(G)$ containing non--trivial $p$--torsion.
Explicitly, we shall have

$\qquad p=2$: $G=G_{2},F_{4},E_{6},E_{7},E_{8}$;

$\qquad p=3$: $G=F_{4},E_{6},E_{7},E_{8}$; and

$\qquad p=5$: $G=E_{8}$.

\noindent In these cases a set of \textsl{generating polynomials} for $\ker
\psi _{p}^{\ast }$ is presented in Propositions 5.5--5.7, and the sets $%
r(G,p)$, $e(G,p)$ and $\{k_{t}\}_{t\in e(G,p)}$ of integers appearing in
Lemma 2.1 are tabulated below, where $e(G,p)$ is given as the subset of $%
r(G,p)$ whose elements are underlined:

\begin{center}
\begin{tabular}{lll}
\hline
$(G,p)$ & $e(G,p)\subset r(G,p)$ & $\{k_{t}\}_{t\in e(G,p)}$ \\ \hline
$(G_{2},2)$ & $\{2,\underline{3}\}$ & $\{2\}$ \\ 
$(F_{4},2)$ & $\{2,\underline{3},8,12\}$ & $\{2\}$ \\ 
$(E_{6},2)$ & $\{2,\underline{3},5,8,9,12\}$ & $\{2\}$ \\ 
$(E_{7},2)$ & $\{2,\underline{3},\underline{5},8,\underline{9},12,14\}$ & $%
\{2,2,2\}$ \\ 
$(E_{8},2)$ & $\{2,\underline{3},\underline{5},8,\underline{9},12,14,%
\underline{15}\}$ & $\{8,4,2,2\}$ \\ 
$(F_{4},3)$ & $\{2,\underline{4},6,8\}$ & $\{3\}$ \\ 
$(E_{6},3)$ & $\{2,\underline{4},5,6,8,9\}$ & $\{3\}$ \\ 
$(E_{7},3)$ & $\{2,\underline{4},6,8,10,14,18\}$ & $\{3\}$ \\ 
$(E_{8},3)$ & $\{2,\underline{4},8,\underline{10},14,18,20,24\}$ & $\{3,3\}$
\\ 
$(E_{8},5)$ & $\{2,\underline{6},8,12,14,18,20,24\}$ & $\{5\}$ \\ \hline
\end{tabular}
.
\end{center}

Combining (2.2) with i) of Lemma 2.1 we get the short exact sequence

\begin{enumerate}
\item[(2.3)] $0\rightarrow H^{even}(BT,G/T;\mathbb{F}_{p})\overset{j}{%
\rightarrow }H^{\ast }(BT;\mathbb{F}_{p})\overset{\psi _{p}^{\ast }}{%
\rightarrow }\frac{H^{\ast }(BT;\mathbb{F}_{p})}{\left\langle \theta
_{s}\right\rangle _{s\in r(G,p)}}\rightarrow 0$.
\end{enumerate}

\noindent It implies that $j$ identifies $H^{even}(BT,G/T;\mathbb{F}_{p})$
with $\ker \psi _{p}^{\ast }=\left\langle \theta _{i}\right\rangle _{i\in
r(G,p)}$. In particular, $\{\theta _{i}\}_{i\in r(G,p)}\subset H^{\ast
}(BT,G/T;\mathbb{F}_{p})$. It follows that, for any pair $\{s,t\}\subset $ $%
r(G,p)$ with $t=s+k(p-1)$, there exists a unique $b_{s,t}\in \mathbb{F}_{p}$
so that a relation of the form

\begin{enumerate}
\item[(2.4)] $\mathcal{P}^{k}(\theta _{s})=b_{s,t}\theta _{t}+\tau _{t}$
with $\tau _{t}\in \left\langle \theta _{s}\right\rangle _{s\in r(G,p),s<t}$
\end{enumerate}

\noindent holds in $H^{\ast }(BT,G/T;\mathbb{F}_{p})$ (resp. in $H^{\ast
}(BT;\mathbb{F}_{p})$ via the injection $j$). Based on the concrete
expression of $\{\theta _{i}\}_{i\in r(G,p)}$ in \S 5.2, the next result is
proved in \S 5.3 by direct computation in the much simpler ring $H^{\ast
}(BT;\mathbb{F}_{p})=\mathbb{F}_{p}[\omega _{1},\cdots ,\omega _{n}]$:

\bigskip

\noindent \textbf{Lemma 2.2. }\textsl{With respect to the degree set }$%
r(G,p) $\textsl{\ of the generating polynomials for }$\ker \psi _{p}^{\ast }$%
\textsl{\ (\S 5.2) given in the table, all non--zero }$b_{s,t}$\textsl{\ in
(2.4) are}

\begin{enumerate}
\item[$p=2:$] $b_{2,3}=1$\textsl{\ for }$G_{2},F_{4},E_{6},E_{7},E_{8}$%
\textsl{;}

$b_{8,12}=1$\textsl{\ for }$F_{4},E_{6},E_{7},E_{8}$\textsl{;}

$b_{3,5}=b_{5,9}=b_{8,9}=1$\textsl{\ for }$E_{6},E_{7},E_{8}$\textsl{;}

$b_{12,14}=1$\textsl{\ for }$E_{7},E_{8}$\textsl{;}

$b_{12,15}=b_{14,15}=1$\textsl{\ for }$E_{8}$\textsl{.}

\item[$p=3:$] $b_{2,4}=1$\textsl{\ for }$F_{4},E_{6},E_{7},E_{8}$\textsl{;}

$b_{6,8}=1$\textsl{\ for }$F_{4},E_{6},E_{7}$\textsl{;}

$b_{4,10}=b_{8,14}=b_{8,10}=1$\textsl{\ for }$E_{7},E_{8}$\textsl{;}

$b_{6,10}=-1$\textsl{\ for }$E_{7}$\textsl{;}

$b_{18,20}=b_{14,20}=b_{18,24}=1$\textsl{\ for }$E_{8}$\textsl{;}

\item[$p=5:$] $b_{k,k+4}=1$\textsl{\ for }$G=E_{8}$\textsl{\ and }$%
k=2,8,14,20$\textsl{.}
\end{enumerate}

\section{$H^{\ast }(G;\mathbb{F}_{p})$ as a module over $\mathcal{A}_{p}$}

In this section we construct $H^{\ast }(G;\mathbb{F}_{p})$ from the
presentation of $H^{\ast }(G/T;\mathbb{F}_{p})$ in ii) of Lemma 2.1, and
relate the $\mathcal{P}^{k}$ action on $H^{\ast }(G;\mathbb{F}_{p})$ to the
values of $b_{s,t}\in \mathbb{F}_{p}$ in (2.4).

The pull back of the universal $T$--bundle $E_{T}\rightarrow BT$ via the
fiber inclusion $\psi $ in (2.1) gives rise to the principle $T$--bundle on $%
G/T$

\begin{enumerate}
\item[(3.1)] $T\rightarrow G\overset{\chi }{\rightarrow }G/T$.
\end{enumerate}

\noindent Since $G/T$ is $1$--connected, the \textsl{Borel} \textsl{%
transgression }

\begin{quote}
$\tau :H^{1}(T;\mathbb{F}_{p})\rightarrow H^{2}(G/T;\mathbb{F}_{p})$
\end{quote}

\noindent defines a basis $\{t_{i}\}_{1\leq i\leq n}$ of $H^{1}(T;\mathbb{F}%
_{p})$ by $\tau (t_{i})=\omega _{i}$. Consequently,

\begin{quote}
$H^{\ast }(T;\mathbb{F}_{p})=\Lambda _{\mathbb{F}_{p}}^{\ast }(t_{1},\ldots
,t_{n})$.
\end{quote}

In the Leray--Serre spectral sequence $\{E_{r}^{\ast ,\ast }(G;\mathbb{F}%
_{p}),d_{r}\}$ of (3.1) one has (\cite{K})

\begin{enumerate}
\item[(3.2)] $E_{2}^{s,t}(G;\mathbb{F}_{p})=H^{s}(G/T)\otimes \Lambda _{%
\mathbb{F}_{p}}^{t}(t_{1},\ldots ,t_{n})$;

\item[(3.3)] the differential $d_{2}:E_{2}^{s,t}(G;\mathbb{F}%
_{p})\rightarrow E_{2}^{s+2,t-1}(G;\mathbb{F}_{p})$ is given by

$\qquad d_{2}(x\otimes t_{k})=x\omega _{k}\otimes 1$, $x\in H^{s}(G/T;%
\mathbb{F}_{p})$, $1\leq k\leq n$.
\end{enumerate}

Over $\mathbb{F}_{p}$ the subring $H^{+}(BT;\mathbb{F}_{p})$ has the
canonical additive basis $\{\omega _{1}^{b_{1}}\cdots \omega
_{n}^{b_{n}}\mid b_{i}\geq 0$, $\sum b_{i}\geq 1\}$. Consider the $\mathbb{F}%
_{p}$--linear map

\begin{enumerate}
\item[(3.4)] $\mathcal{D}:H^{+}(BT;\mathbb{F}_{p})\rightarrow E_{2}^{\ast
,1}(G;\mathbb{F}_{p})=H^{\ast }(G/T;\mathbb{F}_{p})\otimes \Lambda _{\mathbb{%
F}}^{1}$
\end{enumerate}

\noindent by $\mathcal{D}(\omega _{1}^{b_{1}}\cdots \omega
_{n}^{b_{n}})=\omega _{1}^{b_{1}}\cdots \omega _{s}^{b_{s}-1}\cdots \omega
_{n}^{b_{n}}\otimes t_{s}$, where $s\in \{1,\cdots ,n\}$ is the smallest one
such that $b_{s}\geq 1$. Immediate but useful properties of the map $%
\mathcal{D}$ are:

\bigskip

\noindent \textbf{Lemma 3.1.} \textsl{Let }$\beta _{1},\beta _{2}\in
H^{+}(BT;\mathbb{F}_{p})$\textsl{\ and write }$\left[ \theta \right] \in
E_{3}^{s,t}(G;\mathbb{F}_{p})$\textsl{\ for the cohomology class of a }$%
d_{2} $\textsl{--cocycle }$\theta \in E_{2}^{s,t}(G;\mathbb{F}_{p})$\textsl{%
. Then}

\begin{quote}
\textsl{i) }$\mathcal{D}(\ker \psi _{p}^{\ast })\subset \ker d_{2}$\textsl{%
;\quad ii) }$\mathcal{D}(\beta _{1}\beta _{2})-\beta _{1}\mathcal{D}(\beta
_{2})\in \func{Im}d_{2}$\textsl{.}
\end{quote}

\noindent \textsl{In particular, }

\begin{quote}
$\left[ \mathcal{D}(\beta _{1}\beta _{2})\right] =0$\textsl{\ if either }$%
\beta _{1}$\textsl{\ or }$\beta _{2}\in \ker \psi _{p}^{\ast }$\textsl{.}
\end{quote}

\noindent \textbf{Proof. }i) is shown by\textbf{\ }$d_{2}(\mathcal{D}(\theta
))=\theta =0$ in $H^{\ast }(G/T;\mathbb{F}_{p})$ for all $\theta \in \ker
\psi _{p}^{\ast }$. For ii) it suffices to consider the cases where $\beta
_{1},\beta _{2}$ are monomials in $\omega _{1},\cdots ,\omega _{n}$, and the
result comes directly from the definition of $\mathcal{D}$.$\square $

\bigskip

According to i) of Lemma 3.1, the map $\mathcal{D}$ assigns each generating
polynomial $\theta _{s}$ an element

\begin{enumerate}
\item[(3.5)] $\alpha _{2s-1}=:[\mathcal{D}(\theta _{s})]\in E_{3}^{2s-2,1}(G;%
\mathbb{F}_{p})$.
\end{enumerate}

\noindent Since $E_{2}^{s,t}(G;\mathbb{F})=0$ for $s$ odd, one has the%
\textsl{\ canonical }monomorphism

\begin{quote}
$E_{3}^{2k,1}(G;\mathbb{F}_{p})=E_{\infty }^{2k,1}(G;\mathbb{F}_{p})=%
\mathcal{F}^{2k}H^{2k+1}(G;\mathbb{F}_{p})\subset H^{2k+1}(G;\mathbb{F}_{p})$
\end{quote}

\noindent which interprets directly $\alpha _{2s-1}$\ as a cohomology class
of $G$, where $\mathcal{F}$ is the filtration on $H^{\ast }(G;\mathbb{F}%
_{p}) $ induced from $\chi $. Furthermore, by Lemma 3.1 if we write $%
\mathcal{T}$ for the subspace of $H^{\ast }(G;\mathbb{F}_{p})$ spanned by
the set $\left\{ \alpha _{2s-1}\right\} _{s\in r(G,p)}$, the map $\mathcal{D}
$ in (3.4) restricts to a surjection

\begin{enumerate}
\item[(3.6)] $\left[ \mathcal{D}\right] :\ker \psi _{p}^{\ast }=H^{+}(BT,G/T;%
\mathbb{F}_{p})\rightarrow \mathcal{T}\subset H^{odd}(G;\mathbb{F}_{p})$.
\end{enumerate}

Let $\{y_{t}\}_{t\in e(G,p)}$ be the set of $p$--special Schubert classes on 
$G/T$ and put $x_{2t}:=\chi ^{\ast }y_{t}\in H^{2t}(G;\mathbb{F}_{p})$.
Denote by $\Delta (\alpha _{2s-1})_{s\in r(G,p)}$ the $\mathbb{F}_{p}$%
--module in the simple system $\left\{ \alpha _{2s-1}\right\} _{s\in r(G,p)}$
of generators. We formulate $H^{\ast }(G;\mathbb{F}_{p})$ from the
presentation of $H^{\ast }(G/T;\mathbb{F}_{p})$ in ii) of Lemma 2.1, and
relate the $\mathcal{P}^{k}$--action on $H^{\ast }(G;\mathbb{F}_{p})$ to the
coefficients $b_{s,t}\in \mathbb{F}_{p}$ in (2.4).

\bigskip

\noindent \textbf{Lemma 3.2.} \textsl{The inclusion }$\left\{ \alpha
_{2s-1}\right\} _{s\in r(G,p)}$\textsl{, }$\{x_{2t}\}_{t\in e(G,p)}\subset
H^{\ast }(G;\mathbb{F}_{p})$\textsl{\ induces an isomorphism of }$\mathbb{F}%
_{p}$\textsl{--modules}

\begin{quote}
\textsl{i) }$H^{\ast }(G;\mathbb{F}_{p})=\mathbb{F}_{p}[x_{2t}]/\left\langle
x_{2t}^{k_{t}}\right\rangle _{t\in e(G,p)}\otimes \Delta (\alpha
_{2s-1})_{s\in r(G,p)}$\textsl{.}
\end{quote}

\noindent \textsl{Moreover, }$\mathcal{T}$\textsl{\ is an invariant subspace
of all }$P^{k}$\textsl{\ and}

\begin{quote}
\textsl{ii) (2.4) implies that }$P^{k}\alpha _{2s-1}=b_{s,t}\alpha _{2t-1}$%
\textsl{, }$t=s+k(p-1)$.
\end{quote}

\noindent \textbf{Proof.} Assertions i) may be regarded as known, see Ka\v{c}
\cite[Theorem 3]{K} or Ishitoya, Kono, Toda \cite[Theorem 1.1]{IKT}.
However, in the existing literatures there seems no a proof for it. In
particular, the crucial relationship (3.5) between the polynomial $\theta
_{2s}$ and the cohomology class $\alpha _{2s-1}$ was absent. We outline a
proof for it because certain ideas in the process are required by showing
ii).

From ii) of Lemma 2.1 and (3.3) one finds that

\begin{enumerate}
\item[(3.7)] $E_{3}^{\ast ,0}=\func{Im}\chi ^{\ast }=\mathbb{F}%
_{p}[x_{2t}]/\left\langle x_{2t}^{k_{t}}\right\rangle _{t\in e(G,p)}\subset
H^{\ast }(G;\mathbb{F}_{p})$.
\end{enumerate}

\noindent The same method as that used in establishing \cite[Lemma 3.4]{DZ2}
is applicable to show that $E_{3}^{\ast ,1}$\ is spanned by $\left\{ \alpha
_{2s-1}\right\} _{s\in r(G,p)}$ (as a module over $E_{3}^{\ast ,0}$).
Further, since $E_{3}^{\ast ,\ast }$ is generated multiplicatively by $%
E_{3}^{\ast ,0}$ and $E_{3}^{\ast ,1}$ \cite{K, S}, and since

\begin{quote}
$E_{3}^{\dim G-n,n}=E_{2}^{\dim G-n,n}=\mathbb{F}_{p}$
\end{quote}

\noindent (for $E_{2}^{\dim G-n-2,n+1}=E_{2}^{\dim G-n+2,n-1}=0$), we get
from iii) of Lemma 2.1 that

\begin{quote}
$E_{3}^{\ast ,\ast }=\mathbb{F}_{p}[x_{2t}]/\left\langle
x_{2t}^{k_{t}}\right\rangle _{t\in e(G,p)}\otimes \Delta (\alpha
_{2s-1})_{s\in r(G,p)}$.
\end{quote}

\noindent The proof for i) is completed by $E_{3}^{\ast ,\ast }=E_{\infty
}^{\ast ,\ast }=H^{\ast }(G;\mathbb{F}_{p})$, where the first equality comes
from $E_{3}^{\ast ,0}$, $E_{3}^{\ast ,1}\subset H^{\ast }(G;\mathbb{F}_{p})$.

Turning to ii) the short exact sequence (2.3) induces the exact sequence of
complexes

\begin{center}
$0\rightarrow H^{\ast }(BT,G/T;\mathbb{F}_{p})\otimes \Lambda ^{\ast
}\rightarrow H^{\ast }(BT;\mathbb{F}_{p})\otimes \Lambda ^{\ast }\rightarrow 
\mathcal{B}_{2}^{\ast ,\ast }\rightarrow 0$,
\end{center}

\noindent in which

\begin{quote}
$\Lambda ^{\ast }=\Lambda _{\mathbb{F}_{p}}^{\ast }(t_{1},\ldots ,t_{n})$, $%
\mathcal{B}_{2}^{\ast ,\ast }=\frac{H^{\ast }(BT;\mathbb{F}_{p})}{%
\left\langle \theta _{i}\right\rangle _{i\in r(G,p)}}\otimes \Lambda ^{\ast
} $
\end{quote}

\noindent and

\begin{quote}
$H^{\ast }(BT,G/T;\mathbb{F}_{p})\otimes \Lambda ^{\ast }=E_{2}^{\ast ,\ast
}(E_{T},G;\mathbb{F}_{p})$;$\quad $

$H^{\ast }(BT;\mathbb{F}_{p})\otimes \Lambda ^{\ast }=E_{2}^{\ast ,\ast
}(E_{T};\mathbb{F}_{p})$,
\end{quote}

\noindent where $E_{T}$ is the total space of the universal $T$--bundle on $%
BT$. It is clear that $\mathcal{B}_{2}^{\ast ,\ast }$ is a subcomplex of $%
E_{2}^{\ast ,\ast }(G;\mathbb{F}_{p})$ with

\begin{quote}
$\mathcal{B}_{3}^{\ast ,1}=\mathcal{T}$ and $\mathcal{B}_{3}^{\ast ,\ast
}=\Delta (\alpha _{2i-1})_{i\in r(G,p)}\subset H^{\ast }(G;\mathbb{F}_{p})$,
\end{quote}

\noindent Since $E_{3}^{\ast ,\ast }(E_{T};\mathbb{F}_{p})=0$ the connecting
homomorphisms in cohomologies give rise to the isomorphisms

\begin{quote}
$\beta :\mathcal{B}_{3}^{\ast ,1}=\mathcal{T}\rightarrow E_{3}^{\ast
,0}(E_{T},G;\mathbb{F}_{p})$,$\quad \beta ^{\prime }:H^{odd}(G;\mathbb{F}%
_{p})\rightarrow H^{even}(E_{T},G;\mathbb{F}_{p})$
\end{quote}

\noindent that fit in the commutative diagrams

\begin{enumerate}
\item[(3.8)] $%
\begin{array}{ccccc}
0\rightarrow & H^{odd}(G;\mathbb{F}_{p}) & \overset{\beta ^{\prime }}{%
\underset{\cong }{\rightarrow }} & H^{even}(E_{T},G;\mathbb{F}_{p}) & 
\rightarrow 0 \\ 
& \cup &  & \cup \text{ }\kappa &  \\ 
0\rightarrow & \mathcal{T} & \overset{\beta }{\underset{\cong }{\rightarrow }%
} & E_{3}^{even,0}(E_{T},G;\mathbb{F}_{p}) & \rightarrow 0 \\ 
&  & \nwarrow \lbrack \mathcal{D}] & \uparrow \chi ^{\ast } &  \\ 
&  &  & H^{even}(BT,G/T;\mathbb{F}_{p}) & 
\end{array}%
$,
\end{enumerate}

\noindent where the inclusion $\kappa $ identifies $E_{3}^{even,0}(E_{T},G;%
\mathbb{F}_{p})$ with the subring

\begin{quote}
$\func{Im}\chi ^{\ast }[H^{even}(BT,G/T;\mathbb{F}_{p})\rightarrow
H^{even}(E_{T},G;\mathbb{F}_{p})]$.
\end{quote}

\noindent (by a standard property of Leray--Serre spectral sequence). Since $%
[\mathcal{D}]=(\beta ^{\prime })^{-1}\circ \chi ^{\ast }$ by (3.8) and since
both $\beta ^{\prime }$ and $\chi ^{\ast }$ commute with $\mathcal{P}^{k}$,
we obtain ii).$\square $

\bigskip

In the context of \cite[Theorem 1.1]{IKT} the classes $\alpha _{2s-1}\in
H^{odd}(G;\mathbb{F}_{p})$ are called \textsl{transgressive} with \textsl{%
transgressive image} $\theta _{s}$, $s\in r(G,p)$. So it is appropriate for
us to introduced the next definition (in view of i) of Lemma 3.2).

\bigskip

\noindent \textbf{Definition 3.3.} Elements in the set $\left\{ \alpha
_{2s-1}\right\} _{s\in r(G,p)}$ are called $p$\textsl{--transgressive
generators} on $H^{\ast }(G;\mathbb{F}_{p})$.$\square $

\section{Proofs of Theorems 1 and 2}

Assume in this section that $G$ is exceptional with $H^{\ast }(G)$
containing non--trivial $p$--torsion. Let $\left\{ \alpha _{2s-1}=:[\mathcal{%
D}(\theta _{s})]\right\} _{s\in r(G,p)}$ be the set of $p$\textsl{--}%
transgressive generators on $H^{\ast }(G;\mathbb{F}_{p})$ with $\theta _{s}$
being given as those in Proposition 5.5--5.7.

\bigskip

\noindent \textbf{4.1. The Bockstein\ }$\delta _{p}:H^{\ast }(G;\mathbb{F}%
_{p})\rightarrow H^{\ast }(G;\mathbb{F}_{p})$.

Instead of the set $\left\{ \alpha _{2s-1}\right\} _{s\in r(G,p)}$ of $p$%
\textsl{--}transgressive generators on $H^{\ast }(G;\mathbb{F}_{p})$, in 
\cite{DZ2} the ring $H^{\ast }(G;\mathbb{F}_{p})$ was described by the set

\begin{quote}
$\mathcal{O}_{G,\mathbb{F}_{p}}=\left\{ \zeta _{2s-1}\in E_{3}^{2s-2,1}(G;%
\mathbb{F}_{p})\mid s\in r(G,p)\right\} $
\end{quote}

\noindent of $p$\textsl{--primary generators} on $H^{\ast }(G;\mathbb{F}%
_{p}) $ (see \cite[Definition 2.9; Theorems 4--5]{DZ2}). These classes $%
\zeta _{2s-1}$ behave well with respect to the Bockstein\textbf{\ }$\delta
_{p}$ in the sense that (see \cite[Lemma 3.5]{DZ2}):

\begin{enumerate}
\item[(4.1)] $\delta _{p}(\zeta _{2s-1})=\left\{ 
\begin{array}{c}
-x_{2s}\text{ if }s\in e(G,p)\text{;} \\ 
0\text{ if }s\notin e(G,p)\text{.\quad }%
\end{array}%
\right. $
\end{enumerate}

\noindent On the other hand, from the proof of Lemma 3.2 we find that

i) the ring $E_{3}^{\ast ,0}=\mathbb{F}_{p}[x_{2t}]/\left\langle
x_{2t}^{k_{t}}\right\rangle _{t\in e(G,p)}$ has the additive basis

\begin{center}
$\mathcal{C}=\{\tprod\limits_{t\in e(G,p)}x_{2t}^{r(t)}\mid
r:e(G,p)\rightarrow \mathbb{Z}$ is a map with $0\leq r(t)<k_{t}\}$;
\end{center}

ii) $E_{3}^{\ast ,1}(G,\mathbb{F}_{p})$ is a $E_{3}^{\ast ,0}$--module with
basis $\left\{ \alpha _{2s-1}\right\} _{s\in r(G,p)}$.

\noindent As a result each $\xi _{2s-1}\in E_{3}^{2s-2,1}(G;\mathbb{F}_{p})$
can be expressed uniquely by

\begin{enumerate}
\item[(4.2)] $\zeta _{2s-1}=\tsum\limits_{g_{i}\in \mathcal{C}}g_{i}\alpha
_{2t-1}$ with $\deg g_{i}=2(s-t)$.
\end{enumerate}

\noindent Our strategy in establishing the formulae for $\delta _{p}(\alpha
_{2s-1})$ in Theorem 1 is to clarify the expression (4.2), and to apply the
formula (4.1).

\bigskip\ 

\noindent \textbf{Theorem 4.1.} \textsl{We have }$\zeta _{2s-1}=\alpha
_{2s-1}$\textsl{\ with the following exceptions:}

\begin{enumerate}
\item[i)] \textsl{for }$p=2$\textsl{\ and in }$E_{7},E_{8}:$
\end{enumerate}

$\quad \zeta _{15}=\alpha _{15}+x_{6}\alpha _{9}$\textsl{; }$\quad \zeta
_{27}=\alpha _{27}+x_{10}\alpha _{17}$\textsl{\ in }$E_{7},E_{8}$\textsl{,}

$\quad \zeta _{23}=\alpha _{23}+x_{6}\alpha _{17}$\textsl{\ in }$E_{7}$%
\textsl{;}

$\quad \zeta _{23}=\alpha _{23}+x_{6}\alpha _{17}+x_{6}^{3}\alpha _{5}$%
\textsl{; }$\quad \zeta _{29}=\alpha _{29}+x_{6}^{2}\alpha _{17}$\textsl{\
in }$E_{8}$\textsl{.}

\begin{enumerate}
\item[ii)] \textsl{for }$p=3$
\end{enumerate}

$\quad \zeta _{15}=\alpha _{15}-x_{8}\alpha _{7}$\textsl{\ in }$%
F_{4},E_{6},E_{7},E_{8}$\textsl{;}

$\quad \zeta _{35}=\alpha _{35}+x_{8}\alpha _{27}$\textsl{\ in }$E_{7},E_{8}$%
\textsl{;}

$\quad \zeta _{19}=-\alpha _{19};$ $\zeta _{27}=\alpha _{27}+x_{8}\alpha
_{19}$\textsl{;}$\quad $

$\quad \zeta _{39}=\alpha _{39}-x_{20}\alpha _{19}$\textsl{; }$\quad \zeta
_{47}=\alpha _{47}-x_{8}\alpha _{39}$\textsl{\ in }$E_{8}$\textsl{.}

\begin{enumerate}
\item[iii)] \textsl{for }$p=5$\textsl{\ and in }$E_{8}$\textsl{:}

$\zeta _{s}=\left\{ 
\begin{tabular}{l}
$3\alpha _{15}\text{ for }s=15\text{;}$ \\ 
$3\alpha _{23}+2x_{12}\alpha _{11}\text{ for }s=23$; \\ 
$-\alpha _{35}-x_{12}^{2}\alpha _{11}\text{ for }s=35\text{;}$ \\ 
$3\alpha _{47}+x_{12}^{3}\alpha _{11}\text{ for }s=47\text{.}$%
\end{tabular}%
\right. $
\end{enumerate}

\noindent \textbf{Proof.} Let $\Phi _{G,\mathbb{F}_{p}}=\{\gamma
_{s}\}_{s\in r(G,p)}$, $\deg \gamma _{s}=2s$ be the set of $p$--primary
polynomials on $G$ (\cite[Definition 2.7]{DZ2}). In the context of \cite[\S 6%
]{DZ1} each $\gamma _{s}\in \Phi _{G,\mathbb{F}_{p}}$ has been presented as

\begin{enumerate}
\item[(4.3)] $\gamma _{s}=\beta _{s}+\tsum\limits_{r}\beta _{r}y^{r}$ with $%
\beta _{s},\beta _{r}\in \ker \psi _{p}^{\ast }$;
\end{enumerate}

\noindent where

\begin{quote}
i) the sum is over all functions $r:e(G,p)\rightarrow \mathbb{Z}$ with $%
0\leq r(t)<k_{t}$ and $\tsum r(t)>0$,

ii) $y^{r}=\tprod\limits_{t\in e(G,p)}y_{t}^{r(t)}$ with $\{y_{t}\mid t\in
e(G,p)\}$ the set of $p$--special Schubert classes on $G/T$ (see \S 2).
\end{quote}

\noindent Applying the operator $\varphi $ in \cite[(2.7)]{DZ2} to both
sides of (4.3) yields in $E_{3}^{2s-2,1}(G;\mathbb{F}_{p})$ the relation

\begin{enumerate}
\item[(4.4)] $\zeta _{2s-1}=[\varphi (\gamma _{s})]=\mathcal{D}(\beta
_{s})+\tsum x^{r}\mathcal{D}(\beta _{r})$,
\end{enumerate}

\noindent where the first equality comes from the definition of the class $%
\zeta _{2s-1}$ \cite[Definition 2.9]{DZ2}, the second is obtained by
comparing the definitions of $\varphi $ in \cite[(2.7)]{DZ2} with that of $%
\mathcal{D}$ in (2.4), and where $\mathcal{D}(\beta _{s}),$ $\mathcal{D}%
(\beta _{r})\in \mathcal{T}$ by (3.6).

Assume that $\deg \beta _{r}=c$. By i) of Lemma 2.1, $\beta _{s},$ $\beta
_{r}\in \ker \psi _{p}^{\ast }$ implies that

\begin{enumerate}
\item[(4.5)] $\beta _{s}=b_{s}\theta _{s}+\tau _{s}$, $\beta _{r}=\left\{ 
\begin{tabular}{l}
$\tau _{c}\text{ if }c\notin r(G,p)$ \\ 
$b_{r}\theta _{c}+\tau _{c}\text{ if }c\in r(G,p)$%
\end{tabular}%
\right. $,
\end{enumerate}

\noindent where $b_{s},b_{r}\in \mathbb{F}_{p}$, $\tau _{h}\in \left\langle
\theta _{t}\right\rangle _{t\in r(G,p),t<h}$. Consequently

\begin{enumerate}
\item[(4.6)] $\mathcal{D}(\beta _{s})=b_{s}\alpha _{2s-1}$; $\mathcal{D}%
(\beta _{r})=\left\{ 
\begin{tabular}{l}
$0\text{ if }c\notin r(G,p)$ \\ 
$b_{r}\alpha _{2c-1}\text{ if }c\in r(G,p)$%
\end{tabular}%
\right. $.
\end{enumerate}

\noindent by Lemma 3.1. Substituting the expressions (4.6) in (4.4) we get
the desired expression (4.2) of $\zeta _{2s-1}$ in terms of $\alpha _{2s-1}$%
's. This explains the algorithm to obtain the relations in Theorem 4.1.

Finally, we remark that, in the context of \cite{DZ1}, the polynomials $%
\gamma _{s}$ have been explicitly presented in the form of (4.3) (as
examples, see in \cite[(6.2), (6.3)]{DZ1} for the cases $G=E_{7}$ and $p=2,3$%
) and the method computing $b_{s},$ $b_{r}\in \mathbb{F}_{p}$ will be
illustrated by our latter computation for $b_{s,t}\in \mathbb{F}_{p}$
((2.4)) in \S 5.3.$\square $

\bigskip

\noindent \textbf{4.2. Proof of Theorem 1. }The\textbf{\ }presentations of%
\textbf{\ }$H^{\ast }(G;\mathbb{F}_{p})$ in Theorem 1 come from i) of Lemma
3.2, together the degree set $r(G,p)$ given in the table in \S 2. We note
that, in a characteristic $p\neq 2$, the factor $\Delta (\alpha
_{2s-1})_{s\in r(G,p)}$ in Lemma 3.2 can be replaced by the exterior algebra 
$\Lambda (\alpha _{2s-1})_{s\in r(G,p)}$ because odd dimensional cohomology
classes are all square free.

By ii) of Lemma 3.2 results on $\mathcal{P}^{k}(\alpha _{2s-1})$ are
verified by Lemma 2.2.

Combining formula (4.1) with the expressions of $\zeta _{2s-1}$ in Theorem
4.1 yields the formulae of $\delta _{p}(\alpha _{2s-1})$ in Theorem 1. As
examples consider the case $(G,p)=(E_{8},2)$. From Theorem 4.1 we have

\begin{quote}
$\zeta _{2s-1}=\alpha _{2s-1}$, $s=2,3,5,9,$\textsl{\ }

$\zeta _{15}=\alpha _{15}+x_{6}\alpha _{9}$\textsl{; }$\quad \zeta
_{27}=\alpha _{27}+x_{10}\alpha _{17}$\textsl{,}

$\zeta _{23}=\alpha _{23}+x_{6}\alpha _{17}+x_{6}^{3}\alpha _{5}$\textsl{; }$%
\quad \zeta _{29}=\alpha _{29}+x_{6}^{2}\alpha _{17}$\textsl{.}
\end{quote}

\noindent With $e(E_{8},2)=\{3,5,9,15\}\subset
r(E_{8},2)=\{2,3,5,8,9,12,14,15\}$ (see the table in \S 2) we get from (4.1)
that

\begin{quote}
$\delta _{2}\alpha _{3}=0$; $\delta _{2}\alpha _{2s-1}=x_{2s}$; $s=3,5,9$

$\delta _{2}\alpha _{15}+x_{6}\delta _{2}\alpha _{9}=0$; $\delta _{2}\alpha
_{27}+x_{6}\delta _{2}\alpha _{17}=0$;

$\delta _{2}\alpha _{23}+x_{6}\delta _{2}\alpha _{17}+x_{6}^{3}\delta
_{2}\alpha _{5}=0$; $\delta _{2}\alpha _{29}+x_{6}^{2}\delta _{2}\alpha
_{17}=x_{30}$.
\end{quote}

\noindent These justify the formulae for $\delta _{2}\alpha _{2s-1}$ in
(1.1) for the case $G=E_{8}$.$\square $

\bigskip

\noindent \textbf{Proof of Theorem 2. }Let\textbf{\ }$\phi _{p}:H^{\ast }(G;%
\mathbb{F}_{p})\rightarrow H^{\ast }(G;\mathbb{F}_{p})\otimes H^{\ast }(G;%
\mathbb{F}_{p})$ be the reduced co--product induced by the multiplication $%
\mu :G\times G\rightarrow G$. In the notation developed in this paper, Lemma
2.1 in \cite{IKT} can be rephrased as

\begin{enumerate}
\item[(4.7)] $\phi _{p}(\alpha _{2s-1})\in E_{3}^{\ast ,0}\otimes \mathcal{T}
$ ,
\end{enumerate}

\noindent where $E_{3}^{\ast ,0}=\mathbb{F}_{p}[x_{2t}]/\left\langle
x_{2t}^{k_{t}}\right\rangle _{t\in e(G,p)}$. Granted with Theorem 1,
coefficients comparison based on (4.7) suffices to establish Theorem 2. As
examples we show the formulae for $\phi _{p}(\alpha _{2s-1})$ for the cases
of $G=E_{8}$ and $p=2,3$.

With respect to the presentation of $H^{\ast }(E_{8};\mathbb{F}_{2})$ in
(1.1) we can assume that

\begin{quote}
$\phi _{2}(\alpha _{15})=ax_{10}\otimes \alpha _{5}+bx_{6}\otimes \alpha
_{9}+cx_{6}^{2}\otimes \alpha _{3},$ $a,b,c\in \mathbb{F}_{2}$.
\end{quote}

\noindent by (4.7). From $\delta _{2}(\alpha _{15})=x_{6}x_{10}$, $\delta
_{2}(\alpha _{5})=x_{6}$, $\delta _{2}(\alpha _{3})=0$ by Theorem 1 and from
the obvious relation $\delta _{2}\phi _{2}=\phi _{2}\delta _{2}$, one finds
that $a=b=1$. Consequently

\begin{quote}
$\phi _{2}(\alpha _{15})=x_{10}\otimes \alpha _{5}+x_{6}\otimes \alpha
_{9}+cx_{6}^{2}\otimes \alpha _{3}$.
\end{quote}

\noindent Applying $\mathcal{P}^{1}$ to both sides yields the equation

\begin{quote}
$\phi _{2}(\mathcal{P}^{1}\alpha _{15})=x_{6}^{2}\otimes \alpha
_{5}+cx_{6}^{2}\otimes \alpha _{5}$.
\end{quote}

\noindent From $\mathcal{P}^{1}\alpha _{15}=\alpha _{17}=\mathcal{P}^{4}%
\mathcal{P}^{2}\mathcal{P}^{1}\alpha _{3}$ (by (1.1)) and from $\phi
_{2}(\alpha _{3})=0$, we get $c=1$. This verifies the formula for $\phi
_{2}(\alpha _{15})$ in Theorem 2.

Similarly, with respect to the presentation of $H^{\ast }(E_{8};\mathbb{F}%
_{3})$ in (1.2) we can assume by (4.7), for the degree reason, that

\begin{quote}
$\phi _{3}(\alpha _{15})=ax_{8}\otimes \alpha _{7}$

$\phi _{3}(\alpha _{35})=bx_{8}\otimes \alpha _{27}+cx_{8}^{2}\otimes \alpha
_{19}+dx_{8}x_{20}\otimes \alpha _{7}+ex_{20}\otimes \alpha _{15}$,
\end{quote}

\noindent where $a,b,c,d,e\in \mathbb{F}_{3}$. From $\delta _{3}\phi
_{3}=\phi _{3}\delta _{3}$ and from the values of $\delta _{3}(\alpha
_{2s-1})$ in (1.2) we get the equations in $H^{\ast }(E_{8};\mathbb{F}%
_{3})\otimes H^{\ast }(E_{8};\mathbb{F}_{3})$

\begin{quote}
$\qquad -ax_{8}\otimes x_{8}=x_{8}\otimes x_{8};$

$\qquad -bx_{8}\otimes x_{8}x_{20}+cx_{8}^{2}\otimes
x_{20}-dx_{8}x_{20}\otimes x_{8}-ex_{20}\otimes x_{8}^{2}$

$=-x_{8}\otimes x_{8}x_{20}+x_{8}^{2}\otimes x_{20}-x_{8}x_{20}\otimes
x_{8}+x_{20}\otimes x_{8}^{2}$.
\end{quote}

\noindent Coefficients comparison yields that $a=e=-1$; $b=c=$ $d=1$. This
verifies the formula for $\phi _{3}(\alpha _{15})$ and $\phi _{3}(\alpha
_{35})$ in Theorem 2.$\square $

\bigskip

\noindent \textbf{4.3. Applications}: \textbf{the algebra }$H^{\ast }(G;%
\mathbb{F}_{2})$.\textbf{\ }It follows from the proof of Lemma 3.2 that the
set of $2$--transgressive generators on $H^{\ast }(G;\mathbb{F}_{2})$ is
unique. Moreover, one can deduce from (1.1) the next result, that expresses $%
H^{\ast }(G;\mathbb{F}_{2})$ solely by these generators (i.e. without
resorting to the $2$--special Schubert classes on $G/T$).

\bigskip

\noindent \textbf{Corollary 4.2.} \textsl{With respect to the }$2$\textsl{%
--transgressive generators on }$H^{\ast }(G;\mathbb{F}_{2})$\textsl{\ one
has the isomorphisms of algebras}

\begin{quote}
$H^{\ast }(G_{2};\mathbb{F}_{2})=\mathbb{F}_{2}[\alpha _{3}]/\left\langle
\alpha _{3}^{4}\right\rangle \otimes \Lambda _{\mathbb{F}_{2}}(\alpha _{5})$;

$H^{\ast }(F_{4};\mathbb{F}_{2})=\mathbb{F}_{2}[\alpha _{3}]/\left\langle
\alpha _{3}^{4}\right\rangle \otimes \Lambda _{\mathbb{F}_{2}}(\alpha
_{5},\alpha _{15},\alpha _{23})$;

$H^{\ast }(E_{6};\mathbb{F}_{2})=\mathbb{F}_{2}[\alpha _{3}]/\left\langle
\alpha _{3}^{4}\right\rangle \otimes \Lambda _{\mathbb{F}_{2}}(\alpha
_{5},\alpha _{9},\alpha _{15},\alpha _{17},\alpha _{23})$;

$H^{\ast }(E_{7};\mathbb{F}_{2})=\frac{\mathbb{F}_{2}[\alpha _{3},\alpha
_{5},\alpha _{9}]}{\left\langle \alpha _{3}^{4},\alpha _{5}^{4},\alpha
_{9}^{4}\right\rangle }\otimes \Lambda _{\mathbb{F}_{2}}(\alpha _{15},\alpha
_{17},\alpha _{23},\alpha _{27})$;

$H^{\ast }(E_{8};\mathbb{F}_{2})=\frac{\mathbb{F}_{2}[\alpha _{3},\alpha
_{5},\alpha _{9},\alpha _{15}]}{\left\langle \alpha _{3}^{16},\alpha
_{5}^{8},\alpha _{9}^{4},\alpha _{15}^{4}\right\rangle }\otimes \Lambda _{%
\mathbb{F}_{2}}(\alpha _{17},\alpha _{23},\alpha _{27},\alpha _{29})$.
\end{quote}

\noindent \textbf{Proof}. In view of (1.1) it suffices to show that

\begin{enumerate}
\item[(4.8)] $\alpha _{2s-1}^{2}=\left\{ 
\begin{tabular}{l}
$x_{6}$ for $s=2$ and in $G_{2},F_{4},E_{6},E_{7},E_{8}$; \\ 
$x_{4s-2}$ for $s=3,5$ and in $E_{7}$, $E_{8}$; \\ 
$x_{30}+x_{6}^{2}x_{18}$ for $s=8$ and in $E_{8}$,%
\end{tabular}%
\right. $
\end{enumerate}

\noindent and that

\begin{enumerate}
\item[(4.9)] $\alpha _{2s-1}^{2}=0$ for those\textsl{\ }$\alpha _{2s-1}$%
\textsl{\ }belonging to the exterior part.
\end{enumerate}

\noindent These can be deduced directly from $\alpha _{2s-1}^{2}=\delta _{2}%
\mathcal{P}^{s-2}\alpha _{2s-1}$ and (1.1), together with the Adem relation 
\cite{A} and the fact $\mathcal{P}^{s-2}\alpha _{2s-1}\in \mathcal{T}$ by
Lemma 3.2.$\square $

\bigskip

\noindent \textbf{Remark 4.3. }Historically\textbf{\ }the algebras $H^{\ast
}(G;\mathbb{F}_{2})$ (together with the $\mathcal{P}^{k}$ action) were first
obtained by Borel, Araki, Shikata and Thomas \cite{B, Ar1, AS, T}. In terms
of generator and relations, their results agree with those given in
Corollary 4.2. However, in these classical results there is no indication on
the effect of $Sq^{1}$ action on the generators in the \textsl{exterior}
part. The formulae for $\delta _{2}(\alpha _{2s-1})$ in (1.1) implies that
these actions are highly nontrivial:

\begin{quote}
$Sq^{1}(\alpha _{15})=\alpha _{3}^{2}\alpha _{5}^{2}$, $Sq^{1}(\alpha
_{27})=\alpha _{5}^{2}\alpha _{9}^{2}$ in $E_{7},E_{8}$;

$Sq^{1}(\alpha _{23})=\alpha _{3}^{2}\alpha _{9}^{2}$ in $E_{7}$;

$Sq^{1}(\alpha _{23})=\alpha _{3}^{2}\alpha _{9}^{2}+\alpha _{3}^{8}$; $%
Sq^{1}(\alpha _{29})=\alpha _{15}^{2}$ in $E_{8}$.$\square $
\end{quote}

In \cite[Theorem 1]{DZ2} we have presented $H^{\ast }(G;\mathbb{F}_{2})$
additively by the $2$--primary generators $\left\{ \zeta _{2s-1}\right\} $ as

\begin{quote}
$H^{\ast }(G;\mathbb{F}_{2})=\mathbb{F}_{2}[x_{2t}]/\left\langle
x_{2t}^{k_{t}}\right\rangle _{t\in e(G,2)}\otimes \Delta (\zeta
_{2s-1})_{s\in r(G,2)}$.
\end{quote}

\noindent To specify its corresponding algebra structure it suffices to find
the expressions of all the squares $\zeta _{2s-1}^{2}$. This can be done by
i) of Theorem 4.1, (4.8) and (4.9).

\bigskip

\noindent \textbf{Corollary 4.4. }\textsl{With respect to the }$2$\textsl{%
--primary generators on }$H^{\ast }(G;\mathbb{F}_{2})$\textsl{, one has the
isomorphisms of algebras}

\begin{quote}
$H^{\ast }(G_{2};\mathbb{F}_{2})=\mathbb{F}_{2}[x_{6}]/\left\langle
x_{6}^{2}\right\rangle \otimes \Delta _{\mathbb{F}_{2}}(\zeta _{3})\otimes
\Lambda _{\mathbb{F}_{2}}(\zeta _{5})$;

$H^{\ast }(F_{4};\mathbb{F}_{2})=\mathbb{F}_{2}[x_{6}]/\left\langle
x_{6}^{2}\right\rangle \otimes \Delta _{\mathbb{F}_{2}}(\zeta _{3})\otimes
\Lambda _{\mathbb{F}_{2}}(\zeta _{5},\zeta _{15},\zeta _{23})$;

$H^{\ast }(E_{6};\mathbb{F}_{2})=\mathbb{F}_{2}[x_{6}]/\left\langle
x_{6}^{2}\right\rangle \otimes \Delta _{\mathbb{F}_{2}}(\zeta _{3})\otimes
\Lambda _{\mathbb{F}_{2}}(\zeta _{5},\zeta _{9},\zeta _{15},\zeta
_{17},\zeta _{23})$;

$H^{\ast }(E_{7};\mathbb{F}_{2})=\frac{\mathbb{F}_{2}[x_{6},x_{10},x_{18}]}{%
\left\langle x_{6}^{2},x_{10}^{2},x_{18}^{2}\right\rangle }\otimes \Delta _{%
\mathbb{F}_{2}}(\zeta _{3},\zeta _{5},\zeta _{9})\otimes \Lambda _{\mathbb{F}%
_{2}}(\zeta _{15},\zeta _{17},\zeta _{23},\zeta _{27})$;

$H^{\ast }(E_{8};\mathbb{F}_{2})=\frac{\mathbb{F}%
_{2}[x_{6},x_{10},x_{18},x_{30}]}{\left\langle
x_{6}^{8},x_{10}^{4},x_{18}^{2},x_{30}^{2}\right\rangle }\otimes \Delta _{%
\mathbb{F}_{2}}(\zeta _{3},\zeta _{5},\zeta _{9},\zeta _{15},\zeta
_{23})\otimes \Lambda _{\mathbb{F}_{2}}(\zeta _{17},\zeta _{27},\zeta _{29})$%
,
\end{quote}

\noindent \textsl{where}

\begin{quote}
$\zeta _{3}^{2}=x_{6}$ in $G_{2},F_{4},E_{6},E_{7},E_{8}$,

$\zeta _{5}^{2}=x_{10}$, $\zeta _{9}^{2}=x_{18}$ in $E_{7},E_{8}$,

$\zeta _{15}^{2}=x_{30}$; $\zeta _{23}^{2}=x_{6}^{6}x_{10}$ in $E_{8}$.$%
\square $
\end{quote}

\noindent \textbf{Remark 4.5.} Corollary 4.4 was used in \cite[\S 6]{DZ2} to
determine the integral cohomology ring $H^{\ast }(G)$ with respect to the 
\textsl{integral} \textsl{primary generators}.$\square $

\bigskip

\noindent \textbf{Remark 4.6.} In \cite[Theorems 3--5]{DZ2} the algebras $%
H^{\ast }(G;\mathbb{F}_{p})$ were presented by the $p$--primary generators.
Theorems 1, 2 and 4.1 determines the structure of $H^{\ast }(G;\mathbb{F}%
_{p})$ as a Hopf algebra over $\mathcal{A}_{p}$ with respect to these
generators. As examples see \cite[Lemma 3.3]{Du2}, where the result is
applied to determine the near--Hopf ring structure on the integral
cohomology $H^{\ast }(G)$ for all exceptional Lie group $G$.$\square $

\section{Proof of Lemma 2}

In \S 5.1 we obtain formulae for the $\mathcal{P}^{k}$ action on the
universal Chern classes of complex vector bundles. In \S 5.2 we present, for
each exceptional $G$ and prime $p=2,3,5$, a set $\{\theta _{s}\}_{s\in
r(G,p)}$ of generating polynomials for the ideal $\ker \psi _{p}^{\ast }$ in
terms of Chern classes of certain vector bundle on $BT$. With these
preliminaries Lemma 2.2 is establishes in \S 5.3.

\bigskip

\noindent \textbf{5.1. The }$\func{mod}p$\textbf{--Wu formulae. }Let $U(n)$
be the unitary group of rank $n$, and let $BU(n)$ be its classifying space.
It is well known that, for a prime $p$,

\begin{quote}
$H^{\ast }(BU(n),\mathbb{F}_{p})=\mathbb{F}_{p}[c_{1},\ldots ,c_{n}]$
\end{quote}

\noindent where $1+c_{1}+\cdots +c_{n}\in H^{\ast }(BU(n),\mathbb{F}_{p})$
is the \textsl{total Chern class} of the universal complex $n$--bundle $\xi
_{n}$ on $BU(n)$. This implies that each $\mathcal{P}^{k}c_{m}$ can be
written as a polynomial in the $c_{1},\ldots ,c_{n}$, and such expression
may be called the $\func{mod}p$\textsl{--Wu formula} for $\mathcal{P}%
^{k}c_{m}$ \cite{P, Sh}. In the next result we present such formulae for
certain $\mathcal{P}^{k}c_{m}$ that are barely sufficient for a proof of
Lemma 2.2.

\bigskip

\noindent \textbf{Proposition 5.1.} \textsl{The following relations hold in} 
$H^{\ast }(BU(n),\mathbb{F}_{p})$

\begin{enumerate}
\item[i)] $p=2:$

$\mathcal{P}^{r}c_{m}=\tsum\limits_{0\leq t\leq r}\binom{r-m}{t}%
c_{r-t}c_{m+t}$, where $\binom{n}{i}=n(n-1)\cdots (n-i+1)/i!$.

\item[ii)] $p=3:$

$\mathcal{P}^{1}c_{m}=(m+2)c_{m+2}-c_{1}c_{m+1}+(c_{1}^{2}+c_{2})c_{m}$;

$\mathcal{P}%
^{2}c_{m}=c_{2}^{2}c_{m}+c_{1}c_{3}c_{m}-c_{4}c_{m}-c_{1}c_{2}c_{1+m}+(m+1)c_{1}^{2}c_{2+m} 
$

$\qquad +(m-1)c_{2}c_{2+m}-(m+1)c_{1}c_{3+m}+\frac{1}{2}(m^{2}+3m+2)c_{4+m}$;

$\mathcal{P}%
^{3}c_{m}=c_{3}^{2}c_{m}+c_{2}c_{4}c_{m}-c_{1}c_{5}c_{m}+c_{6}c_{m}-c_{2}c_{3}c_{1+m}+c_{5}c_{1+m} 
$

$\qquad
+mc_{2}^{2}c_{2+m}+(1+m)c_{1}c_{3}c_{2+m}-(1+m)c_{4}c_{2+m}-mc_{1}c_{2}c_{3+m} 
$

$\qquad -c_{3}c_{3+m}+\frac{1}{2}(m^{2}+m)c_{1}^{2}c_{4+m}-m^{2}c_{2}c_{4+m}-%
\frac{1}{2}(m^{2}+m)c_{1}c_{5+m}$

$\qquad +\frac{1}{6}(m^{3}+3m^{2}+2m-6)c_{6+m}$

\item[iii)] $p=5:$

$\mathcal{P}^{1}c_{m}=(m+4)c_{m+4}-c_{1}c_{m+3}+(c_{1}^{2}-2c_{2})c_{m+2}$

$\
+(-c_{1}^{3}-2c_{1}c_{2}+2c_{3})c_{m+1}+(c_{1}^{4}+c_{1}^{2}c_{2}+2c_{2}^{2}-c_{1}c_{3}+c_{4})c_{m} 
$.
\end{enumerate}

\noindent \textbf{Proof. }For $p=2$ the expansion of $\mathcal{P}^{r}c_{m}$
comes from the classical Wu--formula \cite{Wu} as $c_{r}\func{mod}2$ is the $%
2r^{th}$ Stiefel--Whitney class of the real reduction of $\xi _{n}$.

For $p>2$ we have the general expansion of $\mathcal{P}^{k}c_{m}$ in terms
of the Schur symmetric functions $s_{\lambda }$ by the formula \cite[ (1.2)]%
{Du1}

\begin{enumerate}
\item[(5.1)] $\mathcal{P}^{k}(c_{m})\equiv \sum\limits_{\lambda
}K_{(1^{m-k},p^{k}),\lambda }^{-1}s_{\lambda }$ $\func{mod}p$,
\end{enumerate}

\noindent where $K_{\mu ,\lambda }^{-1}$ is the \textsl{inverse Kostka number%
} associated to the pair $\{\mu $; $\lambda \}$ of partitions, and where the
sum is over all partitions $\lambda $ of $m+2k(p-1)$. We note in (5.1) that

\begin{enumerate}
\item[(5.2)] for those $(p,k)$ concerned by Proposition 1, \cite[Corollary 2]%
{ER} and \cite[Corollary 5]{Du1} can be applied to evaluate the coefficients 
$K_{(1^{m-k},p^{k}),\lambda }^{-1}$;

\item[(5.3)] each Schur function $s_{\lambda }$ can be expanded as a
polynomial in the $c_{r}$'s by the classical Giambelli formula $s_{\lambda
}=\det (c_{\lambda _{j}^{\prime }+j-i})$ \cite[p.36]{M}, where $\lambda
^{\prime }=(\lambda _{1}^{\prime },\lambda _{2}^{\prime },\ldots )$ is the
partition conjugate to $\lambda $.
\end{enumerate}

\noindent Combining (5.2) and (5.3) one obtains the relations in Proposition
5.1.$\square $

\bigskip

\noindent \textbf{Remark 5.2. }We record below the presentation of (5.1)
from which the relevant inverse Kostka numbers are transparent.\textsl{\ }%
For $p=3$ we have

\begin{quote}
$\mathcal{P}%
^{1}c_{m}=ms_{(1^{m+2})}+s_{(1^{m-1},3)}-s_{(1^{m},2)}-s_{(1^{m-2},2^{2})}$;

$\mathcal{P}%
^{2}c_{m}=s_{(1^{m-2},3^{2})}+(m-1)s_{(1^{m+1},3)}-s_{(1^{m-1},2,3)}-s_{(1^{m-3},2^{2},3)} 
$

$\qquad +\frac{m(m-1)}{2}%
s_{(1^{m+4})}-(m-1)s_{(1^{m+2},2)}-(m-2)s_{(1^{m},2^{2})}+2s_{(1^{m-2},2^{3})} 
$

$\qquad +s_{(1^{m-4},2^{4})}$;

$\mathcal{P}%
^{3}c_{m}=s_{(1^{m-3},3^{3})}+(m-2)s_{(1^{m},3^{2})}-s_{(1^{m-2},2,3^{2})}-s_{(1^{m-4},2^{2},3^{2})} 
$

$\qquad +\frac{(m-1)(m-2)}{2}%
s_{(1^{m+3},3)}-(m-2)s_{(1^{m+1},2,3)}-(m-3)s_{(1^{m-1},2^{2},3)}$

$\qquad +2s_{(1^{m-3},2^{3},3)}+s_{(1^{m-5},2^{4},3)}+\frac{m(m-1)(m-2)}{6}%
s_{(1^{m+6})}-\frac{(m-1)(m-2)}{2}s_{(1^{m+4},2)}$

$\qquad -\frac{(m-2)(m-3)}{2}%
s_{(1^{m+2},2^{2})}+(2m-5)s_{(1^{m},2^{3})}+(m-5)s_{(1^{m-2},2^{4})}$

$\qquad -3s_{(1^{m-4},2^{5})}-s_{(1^{m-6},2^{6})}$
\end{quote}

\noindent For $p=5$ we have

\begin{quote}
$\mathcal{P}%
^{1}c_{m}=ms_{(1^{m+4})}+s_{(1^{m-1},5)}-s_{(1^{m},4)}-s_{(1^{m-2},2,4)}+s_{(1^{m+1},3)} 
$

$\quad
+s_{(1^{m-1},2,3)}+s_{(1^{m-3},2^{2},3)}-s_{(1^{m+3},2)}-s_{(1^{m},2^{2})}-s_{(1^{m-2},2^{3})}-s_{(1^{m-4},2^{4})} 
$.$\square $
\end{quote}

\noindent \textbf{5.2. Generating polynomials for }$\ker \psi _{p}^{\ast }$.

We have seen from Lemma 3.2 the crucial role that the set $\{\theta
_{s}\}_{s\in r(G,p)}$ of polynomials has played in a direct construction of $%
H^{\ast }(G;\mathbb{F}_{p})$ as a Hopf algebra over $\mathcal{A}_{p}$.
Although their full expressions may be seen as lengthy, and as a matter of
fact, their leading terms suffice to establish Lemma 2.2, we choose to
present them in detail for the sake of completeness.

For $n$ indeterminacies $t_{1},\cdots ,t_{n}$ of degree $2$ we set

\begin{enumerate}
\item[(5.4)] $1+e_{1}+\cdots +e_{n}=\dprod\nolimits_{1\leq i\leq n}(1+t_{i})$%
,
\end{enumerate}

\noindent That is, $e_{i}$ is the $i^{th}$ elementary symmetric functions in 
$t_{1},\cdots ,t_{n}$ with degree $2i$. For an exceptional $G$ with rank $n$%
, assume that the set $\{\omega _{i}\}_{1\leq i\leq n}\subset H^{2}(BT)$ of
fundamental weights (see Lemma 2.1) is so ordered as the vertices in the
Dynkin diagram of $G$ in \cite[p.58]{H}. We introduce a set of intermediate
polynomials $c_{k}(G)\in H^{2k}(BT)$ useful in simplifying the our
expression of $\theta _{s}$.

\bigskip

\noindent \textbf{Definition 5.3.} If $G=F_{4}$ we let $c_{k}(F_{4})$, $%
1\leq k\leq 6$, be the polynomial obtained from $e_{k}(t_{1},\cdots ,t_{6})$
in (5.4) by letting

\begin{quote}
$t_{1}=\omega_{4}$;$\quad t_{2}=\omega_{3}-\omega_{4};\quad t_{3}=\omega
_{2}-\omega_{3};$

$t_{4}=\omega_{1}-\omega_{2}+\omega_{3};\quad t_{5}=\omega_{1}-\omega
_{3}+\omega_{4};\quad t_{6}=\omega_{1}-\omega_{4}$.
\end{quote}

\noindent If $G=E_{n}$, $n=6,7,8$, we let $c_{k}(E_{n}),$ $1\leq k\leq n$,
be the polynomial obtained from $e_{k}(t_{1},\cdots ,t_{n})$ in (5.4) by
letting

\begin{quote}
$t_{1}=\omega_{n}$;$\quad t_{2}=\omega_{n-1}-\omega_{n}$;$\quad\cdots$;$%
\quad $

$t_{n-3}=\omega_{4}-\omega_{5};\quad
t_{n-2}=\omega_{3}-\omega_{4}+\omega_{2} $;

$t_{n-1}=\omega _{1}-\omega _{3}+\omega _{2}$;$\quad t_{n}=-\omega
_{1}+\omega _{2}$.$\square $
\end{quote}

\noindent \textbf{Lemma 5.4.} \textsl{The class }$1+c_{1}(F_{4})+\cdots
+c_{6}(F_{4})$\textsl{\ (resp. }$1+c_{1}(E_{n})+\cdots +c_{n}(E_{n})$\textsl{%
, }$n=6,7,8$\textsl{) is the total Chern class of a }$6$\textsl{%
--dimensional (resp. }$n$\textsl{--dimensional) complex bundle }$\xi _{G}$%
\textsl{\ on }$BT$\textsl{.}

\textsl{Moreover, }$c_{1}(G)\in H^{2}(BT)$\textsl{\ can be expressed in
terms of weights as}

\begin{quote}
$c_{1}(G)=\left\{ 
\begin{tabular}{l}
$3\omega _{1}$ for $G=F_{4}$; \\ 
$3\omega _{2}$ for $G=E_{6},E_{7},E_{8}$.%
\end{tabular}%
\right. $
\end{quote}

\noindent \textbf{Proof.} For a $2$--dimensional cohomology class $t\in
H^{2}(BT)$ let $L_{t}$ be the complex line bundle on $BT$ with Euler class $%
t $. Then

\begin{quote}
$\xi _{F_{4}}=\tbigoplus\limits_{1\leq i\leq 6}L_{t_{i}}$ (resp. $\xi
_{E_{n}}=\tbigoplus\limits_{1\leq i\leq n}L_{t_{i}}$, $n=6,7,8$),
\end{quote}

\noindent where $t_{i}$ is the linear form in the weights given in
Definition 5.3.

The expressions of all $c_{r}(G)$ by the special Schubert classes on $G/T$
were deduced in \cite[Lemma 4]{DZ1}, by which the formula for $c_{1}(G)$ is
a special case.$\square $

\bigskip

Let $(G,p)$ be a pair with $H^{\ast }(G)$ containing non--trivial $p$%
--torsion. In Propositions 5.5--5.7 below we present, in accordance to $%
p=2,3,5$, a set $\{\theta _{s}\}_{s\in r(G,p)}$ of generating polynomial for 
$\ker \psi _{p}^{\ast }$ (derived in \cite{DZ1}).

\bigskip

\noindent \textbf{Proposition 5.5.} \textsl{For }$G=G_{2}$\textsl{, }$F_{4}$%
\textsl{\ and }$E_{8}$\textsl{, a set }$\{\theta _{i}\}_{i\in r(G,2)}$%
\textsl{\ of generating polynomials for }$\ker \psi _{2}^{\ast }$\textsl{\
is given by}

\begin{center}
\begin{tabular}{l|l|l|l}
\hline\hline
$\{\theta _{i}\}_{i\in r(G,2)}$ & $G_{2}$ & $F_{4}$ & $E_{8}$ \\ \hline
$\theta _{2}$ & $\omega _{1}^{2}+\omega _{1}\omega _{2}+\omega _{2}^{2}$ & $%
c_{2}$ & $c_{2}$ \\ \hline
$\theta _{3}$ & $\omega _{2}^{3}$ & $c_{3}$ & $c_{3}$ \\ \hline
$\theta _{5}$ &  &  & $c_{5}+\omega _{2}c_{4}$ \\ \hline
$\theta _{8}$ &  & $c_{4}^{2}+\omega _{1}^{2}c_{6}$ & $c_{8}+c_{4}^{2}+%
\omega _{2}^{2}c_{6}+\omega _{2}^{3}c_{5}+\omega _{2}^{8}$ \\ \hline
$\theta _{9}$ &  &  & $\omega _{2}^{2}c_{7}+\omega _{2}c_{8}+\omega
_{2}^{3}c_{6}$ \\ \hline
$\theta _{12}$ &  & $c_{6}^{2}+c_{4}^{3}$ & $c_{6}^{2}+c_{4}^{3}$ \\ \hline
$\theta _{14}$ &  &  & $c_{7}^{2}+c_{4}^{2}c_{6}+\omega _{2}^{2}c_{6}^{2}$
\\ \hline
$\theta _{15}$ &  &  & $c_{7}c_{8}+\ \omega _{2}^{7}c_{8}+\omega
_{2}^{3}c_{4}c_{8}$ \\ \hline
\end{tabular}
\end{center}

\noindent \textsl{and for }$G=E_{6},E_{7}$\textsl{\ by}

\begin{quote}
$\{\theta _{i}\}_{i\in r(E_{6},2)}=$ $\{\theta _{i}\mid
_{c_{7}=c_{8}=0}\}_{i\in r(E_{8},2)\backslash \{14,15\}}$;

$\{\theta _{i}\}_{i\in r(E_{7},2)}=$ $\{\theta _{i}\mid _{c_{8}=0}\}_{i\in
r(E_{8},2)\backslash \{15\}}$.$\square $
\end{quote}

\noindent \textbf{Proposition 5.6.} \textsl{For an exceptional }$G$\textsl{\
with }$G\neq G_{2}$\textsl{, a set }$\{\theta _{i}\}_{i\in r(G,3)}$\textsl{\
of generating polynomials for }$\ker \psi _{3}^{\ast }$\textsl{\ is given by}

\begin{center}
\begin{tabular}{l|l|l|l}
\hline\hline
$\{\theta _{i}\}$ & $F_{4}$ & $E_{6}$ & $E_{7}$ \\ \hline
$\theta _{2}$ & $\omega _{1}^{2}-c_{2}$ & $\omega _{2}^{2}-c_{2}$ & $\omega
_{2}^{2}-c_{2}$ \\ \hline
$\theta _{4}$ & $c_{2}^{2}-c_{4}$ & $c_{2}^{2}-c_{4}$ & $c_{2}^{2}-c_{4}$ \\ 
\hline
$\theta _{5}$ &  & $c_{5}+c_{2}c_{3}$ &  \\ \hline
$\theta _{6}$ & $c_{2}c_{4}-c_{6}$ & $c_{2}c_{4}+c_{3}^{2}-c_{6}$ & $-\omega
_{2}^{3}c_{3}+c_{2}c_{4}-\omega _{2}c_{5}\ +c_{3}^{2}-c_{6}$ \\ \hline
$\theta _{8}$ & $-c_{2}c_{6}$ & $-c_{4}^{2}$ & $-c_{4}^{2}+c_{2}c_{3}^{2}\
-\omega _{2}c_{7}+c_{3}c_{5}$ \\ \hline
$\theta _{9}$ &  & $c_{6}c_{3}$ &  \\ \hline
$\theta _{10}$ &  &  & $-c_{4}c_{3}^{2}+c_{2}c_{3}c_{5}+\
c_{3}c_{7}-c_{5}^{2}$ \\ \hline
$\theta _{14}$ &  &  & $c_{4}c_{5}^{2}+c_{2}c_{5}c_{7}+\ c_{7}^{2}$ \\ \hline
$\theta _{18}$ &  &  & $c_{2}c_{3}^{3}c_{7}+c_{3}^{6}+\
c_{3}^{2}c_{5}c_{7}+c_{3}c_{5}^{3}$ \\ \hline
\end{tabular}

\bigskip

\begin{tabular}{l|l}
\hline\hline
$\{\theta _{i}\}$ & $E_{8}$ \\ \hline
$\theta _{2}$ & $\omega _{2}^{2}-c_{2}$ \\ \hline
$\theta _{4}$ & $c_{2}^{2}-c_{4}$ \\ \hline
$\theta _{8}$ & $-\omega _{2}^{5}c_{3}-\omega _{2}^{3}c_{5}\ -\omega
_{2}^{2}c_{3}^{2}-\omega _{2}^{2}c_{6}-\omega _{2}c_{7}+\ c_{3}c_{5}$ \\ 
\hline
$\theta _{10}$ & $-c_{4}c_{3}^{2}+c_{2}c_{3}c_{5}+\
c_{2}c_{8}+c_{3}c_{7}-c_{5}^{2}$ \\ \hline
$\theta _{14}$ & $c_{4}c_{3}c_{7}+\omega _{2}^{3}c_{3}c_{8}+\
c_{2}c_{3}^{2}c_{6}+c_{2}c_{5}c_{7}\ -\omega
_{2}c_{5}c_{8}-c_{3}^{2}c_{8}+c_{3}c_{5}c_{6}+c_{7}^{2}$ \\ \hline
$\theta _{18}$ & 
\begin{tabular}{l}
$-c_{2}c_{4}^{4}+c_{4}c_{3}^{2}c_{8}+c_{4}c_{6}c_{8}-c_{4}c_{7}^{2}\
-c_{2}c_{3}^{3}c_{7}-c_{2}c_{3}c_{5}c_{8}+\ c_{2}c_{3}c_{6}c_{7}$ \\ 
$-\omega _{2}c_{3}c_{6}c_{8}-c_{3}^{6}\
-c_{3}^{2}c_{6}^{2}-c_{5}c_{6}c_{7}+c_{6}^{3}$%
\end{tabular}
\\ \hline
$\theta _{20}$ & $-c_{2}c_{3}c_{7}c_{8}+\ \omega
_{2}c_{3}c_{8}^{2}+c_{3}^{2}c_{6}c_{8}+c_{5}c_{7}c_{8}$ \\ \hline
$\theta _{24}$ & 
\begin{tabular}{l}
$c_{8}^{3}+c_{2}c_{3}^{2}c_{8}^{2}-\omega
_{2}c_{3}c_{6}^{2}c_{8}+c_{2}c_{3}c_{5}c_{6}c_{8}-c_{3}^{2}c_{5}^{2}c_{8}-%
\omega _{2}c_{3}c_{5}c_{7}c_{8}-c_{3}c_{7}^{3}$ \\ 
$-\omega
_{2}c_{3}c_{6}c_{7}^{2}-c_{2}c_{3}c_{5}c_{7}^{2}+c_{5}^{2}c_{7}^{2}+c_{2}c_{4}^{2}c_{7}^{2}-c_{5}c_{6}^{2}c_{7}-c_{3}^{2}c_{5}c_{6}c_{7}+c_{3}^{4}c_{5}c_{7} 
$ \\ 
$-c_{2}c_{5}^{3}c_{7}-c_{3}^{2}c_{6}^{3}\
+c_{2}c_{4}c_{6}^{3}+c_{3}^{4}c_{6}^{2}$%
\end{tabular}
\\ \hline
\end{tabular}
\end{center}

\noindent \textbf{Proposition 5.7.} \textsl{For }$G=E_{8}$\textsl{, a set of
generating polynomials for }$\ker \psi _{5}^{\ast }$\textsl{\ is given by}

$\theta _{2}=-\omega _{2}^{2}-c_{2}$;

$\theta _{6}=2\omega _{2}^{6}-2\omega _{2}^{3}c_{3}-2\omega _{2}c_{5}\
-2c_{3}^{2}-c_{6}$;

$\theta _{8}=-\omega _{2}^{8}-\omega _{2}^{4}c_{4}\ -2\omega
_{2}^{3}c_{5}-\omega _{2}c_{3}c_{4}-\omega _{2}c_{7}-c_{3}c_{5}\
-c_{4}^{2}-c_{8}$;

$\theta _{12}=-2\omega _{2}^{4}c_{4}^{2}-\omega _{2}^{4}c_{8}+\omega
_{2}^{3}c_{3}^{3}+2\omega _{2}^{3}c_{4}c_{5}-2\omega
_{2}^{2}c_{3}^{2}c_{4}-\omega _{2}^{2}c_{3}c_{7}-2\omega _{2}c_{3}c_{4}^{2}$

$\qquad +c_{3}^{4}-c_{3}c_{4}c_{5}-2c_{5}c_{7}+2c_{6}^{2}$;

$\theta _{14}=-2\omega _{2}^{10}c_{4}+2\omega _{2}^{8}c_{3}^{2}-2\omega
_{2}^{7}c_{7}+\omega _{2}^{5}c_{3}c_{6}-2\omega _{2}^{4}c_{3}c_{7}+2\omega
_{2}^{4}c_{5}^{2}+\omega _{2}^{3}c_{3}^{2}c_{5}$

$\qquad +\omega _{2}^{3}c_{4}c_{7}+\omega _{2}c_{3}c_{4}c_{6}-\omega
_{2}c_{4}^{2}c_{5}+\omega _{2}c_{5}c_{8}-2\omega
_{2}c_{6}c_{7}+c_{3}^{2}c_{4}^{2}-c_{3}^{2}c_{8}$

$\qquad +2c_{3}c_{4}c_{7}+c_{4}^{2}c_{6}+c_{4}c_{5}^{2}+c_{7}^{2}$;

$\theta _{18}=-2\omega _{2}^{8}c_{5}^{2}+2\omega
_{2}^{7}c_{3}^{2}c_{5}-2\omega _{2}^{6}c_{3}^{2}c_{6}+\omega
_{2}^{6}c_{3}c_{4}c_{5}+2\omega _{2}^{5}c_{3}^{2}c_{7}+2\omega
_{2}^{4}c_{3}^{2}c_{8}$

$\qquad +\omega _{2}^{4}c_{4}c_{5}^{2}+2\omega _{2}^{3}c_{3}c_{4}^{3}-\omega
_{2}^{3}c_{3}c_{5}c_{7}+2\omega _{2}^{3}c_{4}^{2}c_{7}-2\omega
_{2}^{3}c_{5}^{3}\ -\omega _{2}^{2}c_{3}^{4}c_{4}-2\omega
_{2}^{2}c_{3}^{3}c_{7}$

$\qquad +\omega _{2}^{2}c_{3}c_{4}^{2}c_{5}+2\omega _{2}^{2}c_{4}^{4}-\omega
_{2}^{2}c_{4}^{2}c_{8}-\omega _{2}c_{3}^{4}c_{5}-2\omega
_{2}c_{3}c_{7}^{2}+\omega _{2}c_{4}^{3}c_{5}-2\omega _{2}c_{4}c_{5}c_{8}$

$\qquad +\omega
_{2}c_{5}^{2}c_{7}-c_{3}^{2}c_{4}c_{8}+c_{3}^{2}c_{5}c_{7}-2c_{3}c_{4}^{2}c_{7}+2c_{3}c_{4}c_{5}c_{6}-c_{3}c_{5}^{3}-2c_{3}c_{7}c_{8}+c_{4}c_{7}^{2} 
$;

$\theta _{20}=-\omega _{2}^{17}c_{3}-\omega _{2}^{13}c_{7}+2\omega
_{2}^{12}c_{4}^{2}+2\omega _{2}^{12}c_{8}+2\omega _{2}^{11}c_{3}c_{6}+\omega
_{2}^{10}c_{3}^{2}c_{4}-\omega _{2}^{9}c_{4}c_{7}$

$\qquad +2\omega _{2}^{8}c_{4}^{3}-\omega _{2}^{7}c_{3}c_{5}^{2}-\omega
_{2}^{6}c_{3}^{3}c_{5}-\omega _{2}^{6}c_{3}^{2}c_{8}+\omega
_{2}^{6}c_{4}c_{5}^{2}-2\omega _{2}^{5}c_{3}^{5}+\omega
_{2}^{5}c_{3}c_{4}^{3}$

$\qquad +\ \omega _{2}^{5}c_{4}^{2}c_{7}+2\omega _{2}^{5}c_{5}^{3}-\omega
_{2}^{4}c_{3}^{4}c_{4}-2\omega _{2}^{4}c_{3}c_{4}^{2}c_{5}-2\omega
_{2}^{4}c_{4}c_{5}c_{7}+\omega _{2}^{3}c_{3}^{4}c_{5}$

$\qquad -2\omega _{2}^{3}c_{3}^{2}c_{4}c_{7}-\omega
_{2}^{3}c_{3}c_{4}c_{5}^{2}+\omega _{2}^{2}c_{3}^{6}+2\omega
_{2}^{2}c_{3}^{2}c_{4}^{3}-\omega _{2}^{2}c_{3}^{2}c_{5}c_{7}-2\omega
_{2}c_{3}^{5}c_{4}$

$\qquad +2\omega _{2}c_{3}^{3}c_{5}^{2}+2\omega
_{2}c_{3}^{2}c_{6}c_{7}+\omega
_{2}c_{4}c_{5}^{3}+2c_{3}^{4}c_{8}+c_{3}^{3}c_{4}c_{7}+c_{3}^{2}c_{7}^{2}+2c_{3}c_{4}^{3}c_{5} 
$

$\qquad +2c_{4}^{5}+c_{4}^{3}c_{8}-2c_{5}^{4}$;

$\theta _{24}=-\omega _{2}^{16}c_{8}-\omega _{2}^{13}c_{3}c_{8}-2\omega
_{2}^{9}c_{3}c_{4}c_{8}+2\omega _{2}^{7}c_{4}c_{5}c_{8}+\omega
_{2}^{6}c_{4}c_{6}c_{8}-2\omega _{2}^{6}c_{5}^{2}c_{8}$

$\qquad +2\omega _{2}^{5}c_{3}c_{8}^{2}+\omega
_{2}^{5}c_{4}c_{7}c_{8}-\omega _{2}^{5}c_{5}c_{6}c_{8}+2\omega
_{2}^{4}c_{4}c_{8}^{2}-\omega _{2}^{4}c_{5}c_{7}c_{8}+\omega
_{2}^{3}c_{3}^{3}c_{4}c_{8}\ $

$\qquad -2\omega _{2}^{3}c_{3}^{2}c_{7}c_{8}+\omega
_{2}^{3}c_{3}c_{4}c_{6}c_{8}-2\omega _{2}^{3}c_{3}c_{5}^{2}c_{8}+\omega
_{2}^{3}c_{6}c_{7}c_{8}+\omega _{2}^{2}c_{4}c_{5}^{2}c_{8}-\omega
_{2}^{2}c_{6}c_{8}^{2}$

$\qquad -2\omega _{2}c_{3}c_{4}c_{8}^{2}-\omega
_{2}c_{4}c_{5}c_{6}c_{8}-2\omega
_{2}c_{7}c_{8}^{2}+c_{3}^{4}c_{4}c_{8}+2c_{3}c_{5}c_{8}^{2}+c_{3}c_{6}c_{7}c_{8} 
$

$\qquad -2c_{5}^{2}c_{6}c_{8}$.$\square $

\bigskip

\noindent \textbf{5.3.} \textbf{Proof of Lemma 2.2. }Granted with the
explicit expressions of the set $\{\theta _{s}\}_{s\in r(G,p)}$ of
generating polynomials for $\ker \psi _{p}^{\ast }$ in Propositions 5.5,
5.6, 5.7 and the $\func{mod}p$ Wu--formulae in Proposition 5.1, we prove
Lemma 2.2.

For $(G,p)=(G_{2},2)$, Lemma 2.2 is directly shown by the computation (see
in Proposition 5.5 for the expressions of $\theta _{2}$, $\theta _{3}$ in $%
G_{2}$)

\begin{quote}
$\mathcal{P}^{1}\theta _{2}=\mathcal{P}^{1}(\omega _{1}^{2}+\omega
_{1}\omega _{2}+\omega _{2}^{2})=\omega _{1}^{2}\omega _{2}+\omega
_{1}\omega _{2}^{2}=\theta _{3}+\omega _{1}\theta _{2}$.
\end{quote}

\noindent So we can assume from now on that $G\neq G_{2}$.

Let $\mathbb{F}_{p}[G]$ be the subring of $H^{\ast }(BT;\mathbb{F}_{p})$
generated by $c_{i}=c_{i}(G)$, the weight $\omega _{r}$ with $r=1$ for $%
F_{4} $, and $r=2$ for $E_{6},E_{7},E_{8}$. Then $\left\{ \theta
_{i}\right\} _{i\in r(G,p)}\subset \mathbb{F}_{p}[G]$ by Propositions
5.5--5.7. Since $c_{r}(G)$'s are the $\func{mod}p$ reduction of the Chern
classes of a vector bundle on $BT$, the Wu--formulae in Proposition 5.1,
together with the Cartan--formula [SE], are applicable to express $\mathcal{P%
}^{k}\theta _{r}$ as an element in $\mathbb{F}_{p}[G]$. It remains for us to
sort out the number $b_{s,t}\in \mathbb{F}_{p}$ in the equation (2.4).

The expressions of $\mathcal{P}^{k}\theta _{r}\in \mathbb{F}_{p}[G]$ may
appear lengthy (in particular, this happens when $G=E_{8}$ and $p=3$ and $5$%
). However, we have two practical methods implementing $b_{s,t}\in \mathbb{F}%
_{p}$. The first utilizes \textsl{Mathematica}, while the second lifts the
computation to an appropriate $S^{1}$--bundle on $BT$ at where, $\theta _{r}$
and $\mathcal{P}^{k}c_{m}$ admit much simpler expressions.

\bigskip

\noindent \textbf{Proof of Lemma 2.2 (Method I).} Based on certain build-in
functions of \textsl{Mathematica }a\textbf{\ }procedure to compute $b_{s,t}$
in (2.4) is given as follows.

For an $i\in r(G,p)$ denote by $\mathcal{G}_{i}(G,p)\subset \mathbb{F}%
_{p}[G] $ a Gr\"{o}bner basis of the ideal generated by the subset $\left\{
\theta _{j}\right\} _{j\in r(G,p),j<i}$. Let $\{s,t\}\subset r(G,p)$ be a
pair with $t=s+k(p-1)$.

\begin{quote}
\textsl{Step 1.} Call GroebnerBasis[ , ] to compute $\mathcal{G}_{t}(G,p)$;

\textsl{Step 2.} Call PolynomialReduce[ , , ] to compute the residue $h_{a}$
of $P^{k}\theta _{s}-a\theta _{t}$\ module $\mathcal{G}_{t}(G,p)$ with $a\in 
\mathbb{F}_{p}$ an indeterminacy;

\textsl{Step 3.} {Take $b_{s,t}=a_{0}$ with }$a_{0}$ the solution to the
equation $h_{a}=0${.}
\end{quote}

\noindent We note that in step 2, the residue $h_{a}$ obtained is always
linear in $a$.$\square $

\bigskip

To demonstrate the second method a few notations are required. Let $\kappa
:S(BT)\rightarrow BT$ be the oriented $S^{1}$--bundle on $BT$ with Euler
class $\omega _{r}\in H^{2}(BT)$, where $r=1$ for $F_{4}$ and $r=2$ for $%
E_{6},E_{7},E_{8}$. Then we have

\begin{quote}
$H^{\ast }(S(BT);\mathbb{F}_{p})=H^{\ast }(BT;\mathbb{F}_{p})\mid _{\omega
_{r}=0}$
\end{quote}

\noindent and the induced ring map $\kappa ^{\ast }$ on cohomology is given
simply by $\kappa ^{\ast }\theta =\theta \mid _{\omega _{r}=0}$.

\bigskip

\noindent \textbf{Example 5.8. }Let $\left\{ \theta _{i}\right\} _{i\in
r(G,p)}$ be the set of generating polynomials for $\ker \psi _{p}^{\ast }$.
Then $\kappa ^{\ast }\theta _{i}$ has simpler expression than that of $%
\theta _{i}$. As an example consider the case $(G,p)=(E_{8},5)$. We get from
Proposition 5.7 that

\begin{quote}
$\kappa ^{\ast }\theta _{2}=-c_{2}$

$\kappa ^{\ast }\theta _{6}=-c_{6}-2c_{3}^{2}$;

$\kappa ^{\ast }\theta _{8}=-c_{8}-c_{3}c_{5}-c_{4}^{2}$;

$\kappa ^{\ast }\theta
_{12}=-2c_{5}c_{7}+2c_{6}^{2}-c_{3}c_{4}c_{5}+c_{3}^{4}$;

$\kappa ^{\ast }\theta _{14}=-c_{3}^{2}c_{8}+c_{7}^{2}\
+2c_{3}c_{4}c_{7}+c_{4}^{2}c_{6}+c_{4}c_{5}^{2}+c_{3}^{2}c_{4}^{2}$;

$\kappa ^{\ast }\theta
_{18}=-2c_{3}c_{7}c_{8}-c_{3}^{2}c_{4}c_{8}+c_{4}c_{7}^{2}+c_{3}^{2}c_{5}c_{7}-2c_{3}c_{4}^{2}c_{7}+2c_{3}c_{4}c_{5}c_{6}-c_{3}c_{5}^{3} 
$;

$\kappa ^{\ast }\theta
_{20}=c_{4}^{3}c_{8}+2c_{3}^{4}c_{8}+c_{3}^{2}c_{7}^{2}+c_{3}^{3}c_{4}c_{7}-2c_{5}^{4}+2c_{3}c_{4}^{3}c_{5}+2c_{4}^{5} 
$;

$\kappa ^{\ast }\theta _{24}=\
2c_{3}c_{5}c_{8}^{2}+c_{3}c_{6}c_{7}c_{8}-2c_{5}^{2}c_{6}c_{8}+c_{3}^{4}c_{4}c_{8} 
$.
\end{quote}

\noindent Moreover, on the subring $\kappa ^{\ast }\mathbb{F}_{5}[E_{8}]=%
\mathbb{F}_{5}[c_{2},\cdots ,c_{8}]$, one has

\begin{quote}
$\mathcal{P}%
^{1}c_{m}=(m+4)c_{m+4}-2c_{2}c_{m+2}+2c_{3}c_{m+1}+(2c_{2}^{2}+c_{4})c_{m}$
\end{quote}

\noindent by Proposition 5.1, where we have reserved $c_{r}$ for $\kappa
^{\ast }c_{r}$, and where $\kappa ^{\ast }c_{1}=0$ by Lemma 5.4.$\square $

\bigskip

The second proof of Lemma 2.2 may appear elaborate, but is useful in
confirming the results obtained from the first one, and may be free of
computer.

\bigskip

\noindent \textbf{Proof of Lemma 2.2 (Method II).} The proof is divided into
two cases in accordance with $\kappa ^{\ast }\theta _{t}=0$ and $\kappa
^{\ast }\theta _{t}\neq 0$.

\noindent \textbf{Case 1.} $\kappa ^{\ast }\theta _{t}=0$. This happens only
when $p=2,t=9$ and $G=E_{6},E_{7},E_{8}$ by Propositions 5.5--5.7. Direct
computation shows that

\begin{quote}
$\mathcal{P}^{1}\theta _{8}=\theta _{9}+\omega _{2}^{4}\theta _{5}$

$\mathcal{P}^{4}\theta _{5}=\theta _{9}+c_{4}\theta _{5}+(\omega
_{2}^{2}c_{4}+c_{6})\theta _{3}+(\omega _{2}^{2}c_{5}+c_{7})\theta _{2}$.
\end{quote}

\noindent These verify the assertions $b_{5,9}=b_{8,9}=1$ in Lemma 2.2.

\noindent \textbf{Case 2. }$\kappa ^{\ast }\theta _{t}\neq 0$. Applying $%
\kappa ^{\ast }$ to the relation (2.4) we get in $H^{\ast }(S(BT);\mathbb{F}%
_{p})$ that

\begin{enumerate}
\item[(5.5)] $\mathcal{P}^{k}(\kappa ^{\ast }\theta _{s})=b_{s,t}\kappa
^{\ast }\theta _{t}+\tau _{t}$ with $\tau _{t}\in \left\langle \kappa ^{\ast
}\theta _{s}\right\rangle _{s\in r(G,p),s<t}$.
\end{enumerate}

\noindent Computation for the case $(G,p)=(E_{8},5)$ is typical enough of
the remaining cases. Carrying on the discussion in the Example 5.8 we find
that

\begin{quote}
$\mathcal{P}^{1}\kappa ^{\ast }\theta _{2}=\kappa ^{\ast }\theta
_{6}+(c_{4}-2c_{2}^{2})\kappa ^{\ast }\theta _{2}$;

$\mathcal{P}^{1}\kappa ^{\ast }\theta _{8}=\kappa ^{\ast }\theta
_{12}+(-c_{2}^{2}+2c_{4})\kappa ^{\ast }\theta
_{8}+(-2c_{3}^{2}+2c_{6})\kappa ^{\ast }\theta _{6}$

$\qquad +(2c_{2}c_{8}+2c_{3}c_{7}-c_{4}c_{6}+2c_{5}^{2})\kappa ^{\ast
}\theta _{2}$

$\mathcal{P}^{1}\kappa ^{\ast }\theta _{14}=\kappa ^{\ast }\theta
_{18}-(c_{2}^{2}c_{3}^{2}+c_{2}c_{8}+c_{3}^{2}c_{4}+\
2c_{3}c_{7}-c_{4}c_{6}+2c_{5}^{2})\kappa ^{\ast }\theta _{8}$

$\qquad +2c_{4}^{3}\kappa ^{\ast }\theta
_{6}-(c_{2}c_{3}^{3}c_{5}-c_{2}c_{3}^{2}c_{4}^{2}+\
2c_{2}c_{3}c_{4}c_{7}+c_{2}c_{4}^{2}c_{6}+c_{2}c_{4}c_{5}^{2}\
-c_{2}c_{7}^{2}$

$\qquad +c_{3}^{2}c_{4}c_{6}+c_{3}c_{4}^{2}c_{5}+\
c_{3}c_{6}c_{7}-c_{4}^{2}c_{8}+2c_{4}c_{5}c_{7}+c_{4}c_{6}^{2}\
-2c_{5}^{2}c_{6}+c_{8}^{2})\kappa ^{\ast }\theta _{2}$

$\mathcal{P}^{1}\kappa ^{\ast }\theta _{20}=\kappa ^{\ast }\theta
_{24}+c_{6}\kappa ^{\ast }\theta
_{18}+(c_{3}^{2}c_{4}-c_{3}c_{7}-c_{4}c_{6}+2c_{5}^{2})\kappa ^{\ast }\theta
_{14}$

$\qquad -(-c_{3}c_{4}c_{5}+c_{4}^{3}-2c_{4}c_{8}+c_{5}c_{7})\kappa ^{\ast
}\theta _{12}$

$\qquad -(c_{2}c_{4}^{2}c_{6}+2c_{3}c_{5}c_{8}+c_{4}^{2}c_{8}\
-c_{4}c_{5}c_{7}+c_{4}c_{6}^{2})\kappa ^{\ast }\theta _{8}$

$\qquad -(c_{2}^{2}c_{7}^{2}+c_{2}c_{3}c_{6}c_{7}+\
2c_{3}^{3}c_{4}c_{5}+c_{3}^{2}c_{4}^{3}+2c_{3}^{2}c_{5}c_{7}\
-c_{3}c_{4}c_{5}c_{6}-c_{3}c_{7}c_{8}$

$\qquad +c_{4}^{2}c_{5}^{2}\ -2c_{5}^{2}c_{8}+2c_{5}c_{6}c_{7})\kappa ^{\ast
}\theta _{6}$

$\qquad -(-2c_{2}c_{4}^{3}c_{8}-c_{2}c_{5}^{4}+c_{2}c_{6}c_{7}^{2}\
-c_{3}^{3}c_{5}c_{8}-c_{3}^{2}c_{4}c_{5}c_{7}+c_{3}c_{4}^{3}c_{7}\
-c_{3}c_{4}^{2}c_{5}c_{6}$

$\qquad +c_{3}c_{5}c_{7}^{2}+c_{3}c_{6}^{2}c_{7}+\
c_{4}^{4}c_{6}+c_{4}^{3}c_{5}^{2}+c_{5}^{3}c_{7})\kappa ^{\ast }\theta _{2}$
\end{quote}

\noindent These imply that $b_{s,s+4}=1$ for $s=2,8,14,20$ by (5.5).$\square 
$

\end{document}